\documentclass[fleqn,11pt,a4paper]{article}
\usepackage{amssymb,latexsym,color,paralist,enumerate}
\usepackage{calrsfs}
\newcommand{\cS}{{\cal S}}

\newcommand{\cQ}{{\cal Q}}

\newcommand{\ovV}{\overline{V}}
\newcommand{\PG}{\mathrm{PG}}
\newcommand{\mrR}{{\mathrm{R}}}
\newcommand{\KK}{\mathbb{K}}

\newcommand{\FF}{\mathbb{F}}

\newcommand{\PP}{\mathbb{P}}

\newcommand{\tV}{\widetilde{V}}

\newcommand{\ve}{\varepsilon}
\newcommand{\vte}{\tilde{\varepsilon}}
\newcommand{\rk}{\mathrm{rank}}
\newcommand{\Aut}{\mathrm{Aut}}

\newcommand{\chr}{\mathrm{char}}
\newcommand{\St}{\mathrm{St}}
\newcommand{\cod}{\mathrm{cod}}
 
\newtheorem{prop}{Proposition}[section]
\newtheorem{cor}[prop]{Corollary}
\newtheorem{defin}[prop]{Definition}
\newtheorem{theo}[prop]{Theorem}
\newtheorem{lemma}[prop]{Lemma}
 
\newtheorem{ex}[prop]{Example}

\newtheorem{note}{Remark}

\newtheorem{conj}[prop]{Conjecture}

\def\<{\langle}
\def\>{\rangle}
\title{Regularity in polar spaces of infinite rank} 
\author{Antonio Pasini}
\date{} 
\begin{document}
\maketitle

\begin{abstract}
In this paper we propose a definition of regularity suited for polar spaces of infinite rank and we investigate to which extent properties of regular polar spaces of finite rank can be generalized to polar spaces of infinite rank. 
\end{abstract}      

\section{Introduction} 

This paper is a continuation of \cite{Pinfty}. In \cite{Pinfty} I mainly stressed on differences between polar spaces of infinite rank and those of finite rank, focusing on properties which hold for all polar spaces of finite rank but fail to hold in many polar spaces of infinite rank, thus unwillingly suggesting that the infinite rank case might be too wild for nice theories can be composed for it. In this paper I support the opposite view. Focusing on regularity and related properties, I shall set down pieces of a theory suited for all polar spaces, including those of infinite rank. As expected, the picture we can get in the infinite rank case is not so neat and simple as for polar spaces of finite rank. Nevertheless, it looks nicer than I dared to hope and more interesting than in the finite rank case.    

\subsection{Premise} 

We are not going to recall all basics on polar spaces; we presume the reader knows them. We only fix here some notation and terminology to be used in this introductory section. More information will be offered in Section \ref{Prel}.
 
We use the symbol $\perp$ to denote collinearity, with the convention that every point is collinear with itself. Given a point $x$ of a polar space $\cS$ we denote by $x^\perp$ the set of points of $\cS$ collinear with $x$ and, for a set of points $X$ of $\cS$, we set $X^\perp := \cap_{x\in X}x^\perp$. Two singular subspaces $X$ and $Y$ are said to be {\em opposite} when $X^\perp\cap Y = Y^\perp\cap X = \emptyset$. In particular, two points are opposite if and only if they are non-collinear and two maximal singular subspaces, henceforth called {\em generators}, are opposite if and only if they are disjoint. 

The symbol $\langle .\rangle$ stands for spans, but we use it for spans in polar spaces as well as in projective and vector spaces. Thus, if $X$ is a set of points of a polar space $\cS$ then $\langle X\rangle$ is the subspace of $\cS$ generated by $X$ and, if $\cS$ admits a projective embedding $\ve:\cS\rightarrow\PG(V)$, then $\langle\ve(X)\rangle$ is the projective subspace of $\PG(V)$ spanned by $\ve(X)$.   

In this paper we will often use an informal simplified notation, writing $\{X,Y\}^\perp$ for $(X\cup Y)^\perp$, $\langle X, Y\rangle$ for $\langle X\cup Y\rangle$, $\langle X\cap Y, Z^\perp\rangle$ for $\langle (X\cap Y)\cup Z^\perp\rangle$ and so on. We trust this notation will not confuse the reader.  

We recall that the rank of a non-degenerate polar space $\cS$, henceforth called $\rk(\cS)$, is the least upper bound of the set of the ranks of the generators of $\cS$, the rank of a projective space being its dimension augmented by 1. Following a well established custom, we write $\rk(\cS) = \infty$ and $\rk(\cS) < \infty$ as shortenings of the sentences ``$\rk(\cS)$ is an infinite cardinal number" and ``$\rk(\cS)$ is finite" respectively. We warn that, while when $\rk(\cS) < \infty$ all generators of $\cS$ have the same dimension, when $\rk(\cS) = \infty$ all generators of $\cS$ are infinite dimensional but in general not all of them have the same dimension.   

All polar spaces to be considered in the sequel are non-degenerate and thick-lined of rank at least 2, by assumption. All embeddings are full projective embeddings. 

\subsection{A definition of regularity}

Different but equivalent ways exists to define regularity for pairs of opposite points of a generalized quadrangle. For instance, the following is a rephrasing of the special case $m = 2$ of a definition stated in \cite[6.4.1]{VM} for generalized $m$-gons: a pair $\{a,b\}$ of opposite points of a generalized quadrangle is {\em regular} if, for any point $c$ opposite both $a$ and $b$, if $|\{a, b,c\}^\perp| > 1$ then $c \in \{a,b\}^{\perp\perp}$; equivalently, if $x, y \in \{a,b\}^\perp$ and $x \neq y$ then $\{x,y\}^\perp = \{a, b\}^{\perp\perp}$. The following is also equivalent to this definition: for every line $\ell$, if $\ell\cap\{a,b\}^\perp\neq\emptyset$ then $\ell\cap\{a,b\}^{\perp\perp}\neq\emptyset$. In particular, in the finite case the latter amounts to say that $|\{a,b\}^{\perp\perp}|$ is equal to the number of lines through a point (compare Payne and Thas \cite[1.3]{PT}). 

A generalization of these definitions to polar spaces of arbitrary but finite rank is proposed in \cite{PVMT} and \cite{CCGP}: two opposite points of a polar space $\cS$ of finite rank are said to form a {\em regular} pair if
\begin{itemize}
\item[$(\mrR1)$] $(X\cup Y)^\perp = \{a,b\}^{\perp\perp}$ for any two opposite generators $X$ and $Y$ of $\{a,b\}^\perp$; equivalently, if $c$ is a point opposite both $a$ and $b$ and $\{a,b,c\}^\perp$ contains two opposite generators of $\{a,b\}^\perp$, then $c\in \{a,b\}^{\perp\perp}$.
\end{itemize}
It is proved in \cite[Lemma 5.5]{PVMT} (also \cite[Proposition 5.1]{CCGP}) that $(\mrR1)$ is equivalent to the following:
\begin{itemize}
\item[$(\mrR2)$] if a generator $M$ of $\cS$ contains a generator of the polar space $\{a,b\}^\perp$, then $M\cap\{a,b\}^{\perp\perp}\neq \emptyset$.
\end{itemize} 
Property $(\mrR2)$ still makes sense when $\rk(\cS) = \infty$. In fact $(\mrR2)$ is involved in a characterization of symplectic polar spaces of arbitrary (possibly infinite) rank (see \cite{CCGP}). Property $(\mrR1)$ also makes sense when $\rk(\cS) = \infty$, however it might possibly be vacuous in certain cases. Indeed it is still an open problem whether polar spaces of infinite rank exist which admit no pair of opposite generators. If $\{a,b\}^\perp$ admits no pair of opposite generators then $(\mrR1)$ is (trivially true but) vacuous for $\{a,b\}$. This is not the unique problem we face with $(\mrR1)$ when $rk(\cS) = \infty$; the following is another one. It is easy to see that  $(\mrR2)$ implies $(\mrR1)$, even if $\rk(\cS) = \infty$. On the other hand, if $(\mrR1)$ holds for a pair $\{a,b\}$ of opposite points and every generator of $\{a,b\}^\perp$ admits an opposite in $\{a,b\}^\perp$, then $(\mrR2)$ holds for $\{a,b\}$ if and only if $X^\perp\cap M\neq\emptyset$ for any generator $X$ of $\{a,b\}^\perp$ and every generator $M$ of $\cS$ containing a generator of $\{a,b\}^\perp$ (see the proofs of \cite[Lemma 5.5]{PVMT} and \cite[Proposition 5.1]{CCGP}). So, if every generator of $\{a,b\}^\perp$ admits an opposite in $\{a,b\}^\perp$ and $\cS$ satisfies the following property $(\mathrm{GS})$, then we can prove that if $(\mrR1)$ holds for $\{a,b\}$ then $(\mrR2)$ also holds, otherwise we get stuck.       
\begin{itemize}
\item[$(\mathrm{GS})$] $X^\perp\cap M \neq \emptyset$ for every generator $M$ and every non-maximal singular subspace $X$ of $\cS$.  
\end{itemize} 
It is well known that $(\mathrm{GS})$ holds true when $\rk(\cS) < \infty$ but in general it fails when $\rk(S) = \infty$. This considered, in \cite{CCGP} we looked at $(\mrR2)$ as a possible definition of regularity when $\rk(\cS)  = \infty$, discarding $(\mrR1)$. In support of this proposal, consider the following sharpening of $(\mrR1)$: 
\begin{itemize}
\item[$(\mrR3)$] we have $(X\cup Y)^\perp = \langle X\cap Y, \{a,b\}^{\perp\perp}\rangle$ for any two generators $X$ and $Y$ of $\{a,b\}^\perp$. 
\end{itemize}
Obviously, $(\mrR3)$ implies $(\mrR1)$. In Section \ref{RDC sec} (Theorem \ref{RA}), without assuming that $\rk(\cS) < \infty$, we shall prove the following: 

\begin{theo}\label{1.1}
Properties $(\mrR2)$ and $(\mrR3)$ are equivalent. 
\end{theo}

As previously recalled, properties $(\mrR1)$ and $(\mrR2)$ are equivalent when $\rk(\cS) < \infty$. In this case $(\mrR3)$ and $(\mrR1)$ are equivalent. We believe that when $\rk(\cS) = \infty$ property $(\mrR3)$ is stronger than $(\mrR1)$, although we have no example at hand which shows that this is indeed the case. 

We are now ready to state our definition of regularity for polar spaces of arbitrary, possibly infinite rank. 

\begin{defin}
\em
We say that a pair of opposite points of $\cS$ is {\em regular} if it satisfies property $(\mrR3)$ (equivalently, $(\mrR2)$). If all pairs of opposite points of  $\cS$ are regular then $\cS$ is said to be {\em regular}.
\end{defin}

\subsection{Tight embeddings}\label{tight} 

We say that an embedding $\ve:\cS\rightarrow\PG(V)$ of a polar space $\cS$ is {\em tight} if $\langle\ve(M\cup M')\rangle = \PG(V)$ for at least one pair of (necessarily opposite) generators $M$ and $M'$ of $\cS$. When $\rk(\cS) = n < \infty$ an embedding $\ve:\cS\rightarrow\PG(V)$ is tight if and only if $\dim(V) = 2n$. If this is the case then $\langle\ve(M\cup M')\rangle = \PG(V)$ for every pair pair of opposite generators $M$ and $M'$ of $\cS$. The following is proved in \cite[Theorem 1.2]{PIIGT}:

\begin{prop}\label{1.2}
An embeddable polar space of finite rank is regular if and only if it admits a tight embedding. 
\end{prop}

In particular, when $\cS$ (has finite rank and) admits a unique embedding, as it is the case when $\cS$ is defined over a division ring of characteristic different from $2$, then $\cS$ is regular if and only if it can be spanned by the union of (any) two opposite generators \cite[Corollary 1.3]{PIIGT}. 

As we shall prove in Section \ref{RE sec} (Theorem \ref{R2}), Proposition \ref{1.2} admits the following generalization:

\begin{theo}\label{1.2 bis}
An embeddable polar space $\cS$ of possibly infinite rank admits a tight embedding if and only if the following holds for at least one (equivalently, every) pair $\{a,b\}$ of opposite points of $\cS$:
\begin{itemize}
\item[$(\mrR4)$] the subspace $\{a,b\}^\perp$ contains two singular subspaces $X$ and $Y$ such that $X$ and $Y$ are opposite generators of $\{a,b\}^\perp$, we have
\begin{equation}\label{1.2 eq 1} 
(X\cup Y)^\perp ~ = ~ \{a,b\}^{\perp\perp}
\end{equation} 
and, if $\ve:\cS\rightarrow \PG(V)$ is an embedding of $\cS$, then
\begin{equation}\label{1.2 eq 2} 
\langle\ve(X \cup Y\cup (X\cup Y)^\perp)\rangle ~ = ~  \PG(V). 
\end{equation}
\end{itemize}
\end{theo}  

When $\cS$ is embeddable the group $\Aut(\cS)$ acts transitively on the set of ordered pairs of opposite points. Accordingly, if $(\mrR4)$ holds for at least one pair $\{a,b\}$ of opposite points of $\cS$ then it holds for any such pair. If furthermore $\rk(\cS) < \infty$ then $\Aut(\cS)$ is also transitive on the set of pairs of opposite singular subspaces of any given rank. Hence, if condition (\ref{1.2 eq 1}) holds for at least one choice of opposite points $a$ and $b$ of $\cS$ and opposite generators $X$ and $Y$ of $\{a,b\}^\perp$, then it holds for any such choice. In this case $(\mrR4)$ is equivalent to $\cS$ being regular. 

In contrast, if $\rk(\cS)$ is infinite then in general $\Aut(\cS)$ is intransitive on the set of pairs of opposite generators and, similarly, the stabilizer of $\{a,b\}$ in $\Aut(\cS)$ acts intransitively on the set of pairs of opposite generators of $\{a,b\}^\perp$. So, it can happen that (\ref{1.2 eq 1}) holds for a particular choice of opposite generators $X$ and $Y$ of $\{a,b\}^\perp$ but not for all of them. If this is the case then $\cS$ is not regular. 

Moreover, when $\rk(\cS) < \infty$ then (\ref{1.2 eq 2}) holds for any choice of opposite singular subspaces $X$ and $Y$. So, when $\rk(\cS) < \infty$ condition (\ref{1.2 eq 2}) is redundant. On the other hand, when $\rk(\cS) = \infty$ then in general $\{a,b\}^\perp$ admits opposite generators which do not satisfy (\ref{1.2 eq 2}) (but might possibly satisfy (\ref{1.2 eq 1})). Notably, this also can happen when $\cS$ is regular. 

Suppose that $\cS$ admits a tight embedding $\ve:\cS\rightarrow\PG(V)$ and $\rk(\cS)$ is infinite. Then in general opposite generators $M$ and $N$ also exist such that $\langle\ve(M\cup N)\rangle \subset \PG(V)$. This fact also explains why, when $\rk(\cS) = \infty$, conditions (\ref{1.2 eq 1}) and (\ref{1.2 eq 2}) of $(\mrR4)$ might hold for certain pairs of opposite generators of $\{a,b\}^\perp$ but not for all of them. Indeed, as one can prove, if $X$ and $Y$ are opposite generators of $\{a,b\}^\perp$ then $\ve(X\cup Y)$ spans $\langle\ve(\{a,b\}^\perp)\rangle$ if and only if $X$ and $Y$ satisfy both properties (\ref{1.2 eq 1}) and (\ref{1.2 eq 2}). 

Summarizing the above discussion, when $\rk(\cS)$ is infinite $(\mrR4)$ is not equivalent to $\cS$ being regular. In fact, as we shall see in Section \ref{Ex more}, polar spaces of infinite rank exist which admit a tight embedding (hence they satisfy $(\mrR4)$)  but are not regular.

\begin{conj}
Every regular polar space admits a tight embedding. 
\end{conj}

\subsection{The three-generators property}

Another characterization of regularity has been obtained in \cite{PVMT} by means of the so-called three-generators property. Let $\rk(\cS) < \infty$. When $\rk(\cS)$  is odd we assume that $\cS$ is thick. So, $\cS$ admits triples of mutually opposite generators and every pair of opposite generators belongs to at least one such triple. Let $M, M_1, M_2$ be such a triple. For every subspace $X$ of $M$ we put
\begin{equation}\label{pi12}
 \pi^M_{1,2}(X) ~:= ~  ((X^\perp\cap M_1)^\perp\cap M_2)^\perp\cap M.
\end{equation}
Thus we obtain a duality $\pi^M_{1,2}$ of $M$. The duality $\pi_{2,1}$ is defined in the same way but switching $M_1$ and $M_2$. Clearly, $\pi_{2,1}^M$ is the inverse of $\pi_{1,2}^M$. Hence $\pi_{1,2}^M$ is a polarity if and only if $\pi_{1,2}^M = \pi_{2,1}^M$. 

\begin{defin}\label{1.3 def} 
\em
We say that a pair $\{M_1, M_2\}$ of opposite generators of $\cS$ satisfies the {\em three-generators property} if $\pi_{1,2}^M$ is a polarity for every generator $M$ of $\cS$ opposite to both $M_1$ and $M_2$.  
\end{defin}

The following has been proved in \cite{PVMT}.

\begin{prop}\label{1.3}
A polar space $\cS$ of finite rank is regular if and only if every pair of opposite generators of $\cS$ satisfies the three-generators property.
\end{prop} 
\begin{note}
\em
Definition \ref{1.3 def} also makes sense when $\cS$ is non-tick of odd rank, but it is vacuous in this case. Indeed in this case no triples of mutually opposite generators exist; consequently, every pair of opposite generators trivially satisfies the 3-generators condition. However $\cS$ is regular in this case (by Proposition \ref{1.2} when $\cS$ is embeddable and \cite[Proposition 5.11]{PVMT} when it isn't). Hence, when  $\cS$ is non-thick of odd rank Proposition \ref{1.3} trivially holds true.  
\end{note} 

The first obstacle we meet when looking for a generalization of Proposition \ref{1.3} suited for polar spaces of infinite rank is the fact that, when $\rk(\cS) = \infty$, the mapping $\pi_{1,2}^M$ defined as in (\ref{pi12}) cannot be a duality. Indeed in this case $\dim(M)$ is infinite and infinite dimensional projective spaces admit no dualities. Moreover, as property $(\mathrm{GS})$ might fail to hold in $\cS$, it can happen that for some point $x\in M$ we have $\pi_{1,2}^M(x) = M$. We must change our setting, replacing dualities with something weaker. We do as follows. 

Given a projective space $\PP$, let $\PP^*$ be its dual. Let $P$ and $P^*$ be the point-set of $\PP$ and the set of hyperplanes of $\PP$ respectively.  

\begin{defin}
\em
A {\em partial duality} of $\PP$ is a mapping $\pi:P\rightarrow P^*\cup\{P\}$ such that, for every choice of distinct points $x$ and $y$ of $\PP$, if both $\pi(x)$ and $\pi(y)$ belong to $P^*$ then $\pi(x) \neq \pi(y)$ and $\pi$ induces a bijection from the line of $\PP$ through $x$ and $y$ to the line of $\PP^*$ through $\pi(x)$ and $\pi(y)$.
\end{defin}

Let $\pi$ be a partial duality of $\PP$. The set $\pi^{-1}(P^*) := \{x\in P~|~\pi(x) \neq P\}$ is a (possibly empty) subspace of $\PP$. If $\pi^{-1}(P^*) = P$ then we say that $\pi$ is {\em non-degenerate}. Obviously, $\pi$ is non-degenerate if and only if it is injective. In contrast, when $\pi^{-1}(P^*) = \emptyset$ we say that $\pi$ is {\em trivial}. If for every two distinct points $x, y \in P$ we have $x\in \pi(y)$ if and only if $y\in \pi(x)$, then we say that $\pi$ is {\em reflexive}. If $\pi$ is reflexive then $P\setminus\pi^{-1}(P^*) \subseteq \pi(x)$ for every $x\in \pi^{-1}(P^*)$. However no hyperplane of $\PP$ contains the complement of a proper subspace of $\PP$. Therefore $\pi$ is reflexive only if it is either non-degenerate or trivial. A non-trivial reflexive partial duality is called a {\em polarity}. 

When $\dim(\PP) < \infty$ the dualities of $\PP$ are precisely the (mappings from the poset of subspaces of $\PP$ to the poset of subspaces of $\PP^*$ defined by) non-degenerate partial dualities; a duality is involutory if and olny if it is reflexive.  

Turning back to our polar space $\cS$ but now allowing $\rk(\cS)$ to be infinite, let $M, M_1, M_2$ be three mutualy opposite generators of $\cS$. Equation (\ref{pi12}), whith $X$ ranging in the set of points of $M$ instead of the set of subspaces of $M$, defines a partial duality $\pi_{1,2}^M$ of $M$ (Section \ref{PD}, Lemma \ref{3G0}). As in the finite rank case, $\pi_{1,2}^M$ is a polarity if and only if $\pi_{1,2}^M = \pi_{2,1}^M$ (Section \ref{PD}, Lemma \ref{3G1}). Definition \ref{1.3 def} still makes sense in the present setting. We are not going to repeat it here. However we need one more definition. 

\begin{defin}\label{1.3 def hyp}
\em
A generator $M$ of $\cS$ is {\em hyperbolic} if every hyperplane of $M$ is contained in at most one generator other than $M$.
\end{defin}

As we shall prove in Section \ref{3G sec} (Theorem \ref{3G}), the following holds:

\begin{theo}\label{3G intro}
Let $\cS$ be a polar space of infinite rank, defined over a division ring of characteristic different from $2$. Let $M_1$ and $M_2$ be opposite generators of $\cS$ such that at least one generator of $\cS$ exists which is opposite to both $M_1$ and $M_2$. Suppose moreover that neither $M_1$ nor $M_2$ are hyperbolic. Under these hypoyheses, $\{M_1, M_2\}$ satisfies the three-generators property if and only if $\langle M_1\cup M_2\rangle = \cS$.
\end{theo}    
 
\begin{conj}\label{3G conj} 
The hypothesis that the underlying division ring of $\cS$ has characteristic different from $2$ can be dropped from Theorem \ref{3G intro}, provided that the conclusion of that theorem is rephrased as follows: the pair $\{M_1, M_2\}$ satisfies the three-generators property if and only if $\cS$ admits a (necessarily tight) embedding $\ve:\cS\rightarrow \PG(V)$ such that $\langle \ve(M_1\cup M_2)\rangle = \PG(V)$. 
\end{conj} 
\textbf{Organization of the paper.} In Section \ref{Prel} we state some terminology and recall some basics on projective embeddings of polar spaces. Section \ref{Reg} contains the proofs of Theorems \ref{1.1} and \ref{1.2 bis} and more results on regularity and tight embeddings. Section \ref{3G game} contains the proof of Theorem \ref{3G intro} and one more result in the same vein as Theorem \ref{3G intro}. A number of examples of polar spaces of infinite rank are discussed in Section \ref{Examples}. Most of them are regular.   

\section{Preliminaries}\label{Prel}

As in the Introduction, throughout this section $\cS$ is a non-degenerate thick-lined polar space of rank at least 2, possibly $\rk(\cS) = \infty$. 

\subsection{A survey of elementary properties}

In this subsection we recall a few well known properties of subspaces of $\cS$, to be be freely used in the sequel of this paper. 

We have $\langle X\rangle  \subseteq X^{\perp\perp}$ for every set $X$ of points of $\cS$. Also, if $X\subseteq Y$ then $X^\perp\supseteq Y^\perp$. Consequently, $X^\perp  = X^{\perp\perp\perp}$ for every set $X$ of points of $\cS$. If moreover $X$ is a singular subspace of $\cS$ then $X = X^{\perp\perp}$ if and only if $X$ is the intersection of a family of generators (as it is always the case when $\rk(\cS) < \infty$). Note that, in general, when $\rk(\cS) = \infty$ non-maximal singular subspaces exist which are contained in a unique generator. If $X$ is one of them then $X^{\perp\perp}$ is the unique generator containing $X$.  

Let $X$ be a singular subspace of $\cS$ and $M$ a generator containing $X$. Let $\cod_M(X)$ be the codimension of $X$ in $M$. If $\cod_M(X) < \infty$ then $\cod_N(X) = \cod_M(X)$ for every generator $N$ containing $X$. Indeed the star of $X$ is a (possibly degenerate) polar space of finite rank equal to $\cod_M(X)$. If $\cod_M(X) = 1$ for some (hence every) generator $M$ containing $X$ then we say that $X$ is a {\em sub-generator} of $\cS$. 

Let $a$ and $b$ be two opposite points of $\cS$. Then $\{a,b\}^\perp \cong \St(c)$ for every $c\in \{a,b\}^{\perp\perp}$, where $\St(c)$ (the {\em star} of $c$) is the polar space with the lines and the planes of $\cS$ through $c$ as points and lines respectively. All generators of $\{a,b\}^\perp$ are sub-generators of $\cS$. So, if $N$ and $N'$ are generators of $\{a,b\}^\perp$ then  $M = \langle N, a\rangle$ and $M' = \langle N', b\rangle$ are generators of $\cS$. Note that $M\cap M' = N\cap N'$. Indeed, if $c\in M\cap M'$ then $c\in \{a, b\}^\perp\cap N^\perp\cap N'^\perp = N\cap N'$. Hence $M$ and $M'$ are opposite if and only if $N$ and $N'$ are opposite. Conversely, if $M$ and $M'$ are opposite generators of $\cS$ and $a\in M$, $b\in M'$ are opposite, then $N := M\cap b^\perp$ and $N' := M'\cap a^\perp$ are opposite generators of $\{a,b\}^\perp$. 

Sets as $\{a,b\}^{\perp\perp}$ for $a$ and $b$ opposite points are called {\em hyperbolic lines}. Note that the points of a hyperbolic line $\{a,b\}^{\perp\perp}$ are mutually opposite (in particular, $\{a, b\}^\perp\cap\{a,b\}^{\perp\perp} = \emptyset$) and if $c$ and $d$ are distinct points of $\{a,b\}^{\perp\perp}$ then $\{c,d\}^{\perp\perp} = \{a,b\}^{\perp\perp}$. Note also that $\langle \{a,b\}^\perp, c\rangle = c^\perp$ for every point $c\in \{a,b\}^{\perp\perp}$.  Accordingly, $\langle \{a,b\}^\perp\cup\{a,b\}\rangle$ contains both $a^\perp$ and $b^\perp$. However $a^\perp$ and $b^\perp$ are distinct maximal subspaces of $\cS$. Hence $\langle \{a,b\}^\perp\cup\{a,b\}\rangle = \cS$.    

\subsection{Projective embeddings of point-line geometries} 

In the present subsection we recall some generalities and fix some terminology on projective embeddings of point-line geometries. This will be helpful in the next subsection, where projective embeddings of polar spaces will be discussed.

As in Shult \cite{S1}, a projective embedding of a point-line geometry ${\cal G} = (P, {\cal L})$, henceforth called just {\em embedding} of $\cal G$ for short, is an injective mapping $\ve:P\rightarrow\PP$ from the point-set $P$ of $\cal G$ to the point-set of a projective geometry $\PP$ such that $\ve(P)$ spans $\PP$ and $\ve(\ell) := \{\ve(p)\}_{p\in\ell}$ is a line of $\PP$ for every line $\ell\in{\cal L}$ of $\cal G$. We take the dimension of $\PP$ as the {\em dimension} $\dim(\ve)$ of $\ve$. Note that if $\cal G$ admits skew lines then necessarily $\dim(\PP) > 2$, hence $\PP$ is desarguesian. If $\PP$ is desarguesian and $\KK$ is the underlying division ring of $\PP$ then we say that $\ve$ is {\em defined over} $\KK$. If all embeddings of $\cal G$ are defined over the same division ring $\KK$ we say that $\cal G$ is {\em defined over} $\KK$. 

Given two embeddings $\ve:{\cal G}\rightarrow\PP$ and $\ve':{\cal G}\rightarrow\PP'$ of $\cal G$ defined over the same division ring $\KK$, a {\em morphism} from $\ve$ to $\ve'$ is a morphism of projective geometries $\varphi:\PP\rightarrow\PP'$ such that $\ve' = \varphi\circ\ve$. (See Faure and Fr\"{o}licher \cite[Chaper 6]{FF}) for morphisms of projective geometries.) Note that the condition $\ve' = \varphi\circ\ve$ forces the morphism $\varphi:\PP\rightarrow\PP'$ to be surjective (and, when $\cal G$ is connected, it uniquely determines $\varphi$). If $\varphi$ is also injective then we say that $\varphi$ is an {\em isomorphism} from $\ve$ to $\ve'$ and we write $\ve\cong \ve'$. In general, if a morphism exists from $\ve$ to $\ve'$ we say that $\ve$ {\em covers} $\ve'$ and $\ve'$ is a {\em quotient} of $\ve$. A motivation for this terminology is the following: given a morphism $\varphi:\ve\rightarrow\ve'$, let $K := \mathrm{Ker}(\varphi)$ be the kernel of $\varphi$ (notation and terminology as in \cite{FF}), let $p_K$ be the projection of $\PP$ onto the star $\PP/K$ of $K$ in $\PP$ and put $\ve/K := p_K\circ\ve$. Then $\PP' \cong \PP/K$ and $\ve/K\cong \ve'$. 

Following Shult \cite{S1}, we say that an embedding is {\em relatively universal} if it admits no proper cover. As proved by Ronan \cite{Ron}, every embedding $\ve$ of a geometry $\cal G$ is covered by a relatively universal embedding $\vte$ of $\cal G$, uniquely determined by $\ve$ up to isomorphisms and characterized by the following property: $\vte$ covers all embeddings which cover $\ve$. We call $\vte$ the {\em hull} of $\ve$. Thus, an embedding is relatively universal if and only if it is its own hull. 

An embedding of $\cal G$ is said to be {\em absolutely universal} if it covers all embeddings of $\cal G$. So, $\cal G$ admits the absolutely universal embedding if and only if all embeddings of $\cal G$ admit the same hull; equivalently, up to isomorphism, $\cal G$ admits a unique relatively universal embedding. Clearly, if $\cal G$ admits the absolutely universal embedding and at least one of its embeddings is defined over a given division ring $\KK$, then $\cal G$ itself is defined over $\KK$. 

Finally, we say that an embedding is {\em minimal} if it admts no proper quotients. For instance, tight embeddings of polar spaces, as defined in Section \ref{tight}, are minimal.  

\subsection{Projective embeddings of polar spaces}  

As proved by Tits \cite[chapters 8 and 9]{T} (see also Buekenhout and Cohen \cite[chapters 7-11]{BC} and Cuypers et al. \cite{CJP1}) all polar spaces of rank at least 3 are embeddable but for two families of polar spaces of rank 3. The non-embeddable ones are the line-grassmannians of 3-dimensional projective spaces defined over non-commutative division rings and certain polar spaces with Moufang but non-desarguesian planes, described in \cite[Chapter 9]{T}. 

Let $\cal S$ be an embeddable polar space and, when $\rk(\cS) = 2$, assume that $\cS$ is neither a grid of order at least $4$ nor a generalized quadrangle as in Tits \cite[8.6(II)(a)]{T}. Then $\cS$ admits the absolutely universal embedding as well as a unique minimal embedding (Tits \cite[Chapter 8]{T}; see also Johnson \cite{J1} and \cite{J2} and Cuypers et al. \cite{CJP2}). In the two excluded cases, all embeddings are 3-dimensional, hence they are both relatively universal and minimal. 

In any case, all embeddings of $\cS$ have dimension at least $2\cdot\rk(\cS)-1$ ($\geq 3$ as $\rk(\cS) \geq 2$ by assumption); hence they embed $\cS$ in desarguesian projective spaces. When $\cS$ is not an infinite grid, these projective spaces are defined over the same division ring (so $\cS$ is defined over that ring).

Let $\ve:\cS\rightarrow\PP = \PG(V)$ be an embedding of a polar space $\cS$. The $\ve$-image $\ve(\cS)$ of $\cS$ is a full subgeometry of $\PG(V)$ and there exists a unique quasi-polarity $\pi_\ve$ of $\PG(V)$ (see Buekenhout and Cohen \cite[Definition 7.1.9]{BC} for the definition of quasi-polarities) such that all points of $\ve(\cS)$ are absolute for $\pi_\ve$ and for any two points $x, y \in \cS$ we have $x\perp y$ if and only if $\ve(x)\perp_\ve\ve(y)$, where $\perp_\ve$ is the orthogonality relation associated to $\pi_\ve$ (Buekenhout and Cohen \cite[Chapter 9]{BC}). 

More explictly, let $\cS_{\pi_\ve}$ be the (possibly degenerate) polar space defined by $\pi_\ve$ on $\PG(V)$. Then $\ve(\cS)$ is a subspace of $\cS_{\pi_\ve}$, possibly $\ve(\cS) = \cS_{\pi_\ve}$. Note that, as $\cS$ is non-degenerate by assumption, if $\ve(\cS) = \cS_{\pi_\ve}$ then $\pi_\ve$ is a polarity, namely its radical is trivial. (Recall that the {\em radical} of a quasi-polarity $\pi$ of a projective geometry $\PP$ is the subspace of $\PP$ formed by the points $p\in \PP$ such that $\pi(p) = \PP$.) 

Suppose that $\ve$ is relatively universal and let $\KK$ be its underlying division ring. Then one of the following occurs (Tits \cite[Chaper 8]{T}).  
\begin{itemize}
\item[(1)] $\chr(\KK) \neq 2$, the polarity $\pi_\ve$ is defined by a non-degenerate alternating form and $\ve(\cS) = \cS_{\pi_\ve}$.  
\item[(2)] The quasi-polarity $\pi_\ve$ is defined by the sesquilinearized of a non-degenerate $\sigma$-quadratic form $q:V\rightarrow \KK/\KK_{\sigma,1}$ as defined by Tits \cite[Chapter 8]{T} for an involutory anti-automorphism $\sigma$ of $\KK$ and $\ve(\cS)$ is the polar space $\cS_q$ associated to $q$. In this case, if either $\chr(\KK) \neq 2$ or $\sigma$ acts non-trivially on the center $Z(\KK)$ of $\KK$, then $\pi_\ve$ is a polarity and $\ve(\cS) = \cS_q = \cS_{\pi_\ve}$.  
\end{itemize}
Suppose moreover that $\ve$ is absolutely universal (as when $\rk(\cS) > 2$ or case (1) occurs). In case (1) and in case (2) with $\pi_\ve$ a polarity $\ve$ is the unique embedding of $\cS$. When $\pi_\ve$ is not a polarity we can factorize $\ve$ over any subspace $K$ of the radical $R_\ve$ of $\pi_\ve$, thus obtaining quotients $\ve/K$ of $\ve$. In particular, $\ve/R_\ve$ is the minimum embedding of $\cS$.

\subsubsection{Relations between $\perp$ and $\perp_\ve$}

Let $X$ be a set of points of $\cS$ and $\ve:\cS\rightarrow\PG(V)$ an embedding of $\cS$. Clearly $\langle \ve(X^\perp)\rangle \subseteq \ve(X)^{\perp_\ve}$.  

\begin{prop}\label{perpperp1}
If $X^\perp\not\subseteq X^{\perp\perp}$ then $\ve(X)^{\perp_\ve} = \langle \ve(X^\perp)\rangle$.
\end{prop}
{\bf Proof.} Let $\cal X$ be a projective subspace of $\PG(V)$ and suppose that ${\cal X}\cap\ve(\cS) \not\subseteq {\cal X}^{\perp_\ve}$. Then $\cal X$ is spanned by ${\cal X}\cap\ve(\cS)$. When $\ve$ is relatively universal this claim follows from Tits \cite[\S\S 8.1.6, 8.2.7]{T}. Otherwise, $\ve$ is a quotient of the absolutely universal embedding $\vte:\cS\rightarrow \PG(\widetilde{V})$ of $\cS$, say $\ve = \vte/K$ for a subspace $K$ of the radical $R_{\vte}$ of $\pi_{\vte}$. Let $p_K$ be the projection of $\PG(\widetilde{V})$ onto $\PG(V) = \PG(\widetilde{V})/K$ and put $\widetilde{\cal X} := p_K^{-1}({\cal X})$. Then ${\cal X}\cap\ve(\cS) = p_K(\widetilde{\cal X}\cap\vte(\cS))$. However ${\cal X}\cap\ve(\cS) \not\subseteq {\cal X}^{\perp_\ve}$ by assumption. Hence $\widetilde{\cal X}\cap\vte(\cS) \not\subseteq {\cal X}^{\perp_{\vte}}$. By Tits \cite[\S\S 8.1.6, 8.2.7]{T}, the set $\widetilde{\cal X}\cap\vte(\cS)$ spans $\widetilde{\cal X}$. Hence ${\cal X}\cap\ve(\cS)$ spans $\cal X$. 

Let now $X^\perp\not\subseteq X^{\perp\perp}$ as in the hypotheses of the lemma and put ${\cal X} = \ve(X)^{\perp_\ve}$. Then ${\cal X}\cap\ve(\cS) = \ve(X^\perp)$, which is not contained in ${\cal X}^{\perp_\ve}$ because $X^\perp\not\subseteq X^{\perp\perp}$. By the previous paragraph, ${\cal X}\cap\ve(\cS)$ spans $\cal X$, namely $\ve(X^\perp)$ spans $\ve(X)^{\perp_\ve}$. \hfill $\Box$    

\begin{cor}\label{perpperp2}
If $X^\perp\not\subseteq X^{\perp\perp}\not\subseteq X^\perp$ then
$\langle \ve(X^{\perp\perp})\rangle = \ve(X)^{\perp_\ve\perp_\ve}$.
\end{cor}
{\bf Proof.} Let $X^\perp\not\subseteq X^{\perp\perp}$. Then $\langle \ve(X^\perp)\rangle = \ve(X)^{\perp_\ve}$ by Proposition \ref{perpperp1}. Suppose moreover that $X^{\perp\perp}\not\subseteq X^\perp$. Then, since $X^{\perp} = X^{\perp\perp\perp}$, we have that $\langle \ve(X^{\perp\perp})\rangle  = \ve(X^\perp)^{\perp_\ve}$ by Lemma \ref{perpperp1} with $X^\perp$ in place of $X$. However $\ve(X^\perp)^{\perp_\ve} = \langle \ve(X^\perp)\rangle^{\perp_\ve} = \ve(X)^{\perp_\ve\perp_\ve}$, since $\langle \ve(X^\perp)\rangle = \ve(X)^{\perp_\ve}$. Finally $\langle \ve(X^{\perp\perp})\rangle = \ve(X)^{\perp_\ve\perp_\ve}$, as claimed. \hfill $\Box$  

\begin{cor}\label{perpperp3}
We have $\langle \ve(\{a,b\}^\perp)\rangle = \{\ve(a), \ve(b)\}^{\perp_\ve}$ and  $\langle \ve(\{a,b\}^{\perp\perp})\rangle = \{\ve(a), \ve(b)\}^{\perp_\ve\perp_\ve}$ for any two opposite points $a$ and $b$ of $\cS$.  
\end{cor}
{\bf Proof.} We have $\{a,b\}^\perp\cap\{a,b\}^{\perp\perp} = \{a,b\}^{\perp\perp\perp}\cap\{a,b\}^{\perp\perp} = \emptyset$. The conclusion follows from Proposition \ref{perpperp1} and Corollary \ref{perpperp2}.  \hfill $\Box$ 

\subsubsection{Subspaces of an embeddable polar space}

Given an embedding $\ve:\cS\rightarrow\PG(V)$, a subspace $X$ of $\cS$ {\em arises} from $\ve$ if $X = \ve^{-1}({\cal X})$ for a projective subspace $\cal X$ of $\PG(V)$; equivalently, $\langle \ve(X)\rangle\cap\ve(\cS) = \ve(X)$.

\begin{defin}
\em 
Given a singular subspace $K$ of $\cS$ (possibly $K = \emptyset$) let $\{X_i\}_{i\in I}$ be a family of singular subspaces such that each of them contains $K$ as a hyperplane, $X_i^\perp\cap X_j = K$ for any choice of $i, j \in I$ with $i\neq j$ and $|I| > 1$. Then $X := \cup_{i\in I}X_i$ is a subspace of $\cS$, henceforth called a {\em rosette}. Clearly, $K$ is the radical of $X$. In particular, if $K = \emptyset$ then $X$ is a set of mutually opposite points. 
\end{defin} 

\begin{prop}\label{App1}
Let $\cS$ be embeddable and let $X$ be a subspace of $\cS$. Then $X$ arises from an embedding of $\cS$, except possibly when $X$ a rosette.
\end{prop}
{\bf Proof.} When $\rk(\cS) < \infty$ the above is just the main result of \cite{CGP}. Suppose that $\cS$ has infinite rank. Then $\cS$ admits the universal embedding, say $\ve:\cS\rightarrow \PG(V)$. We shall prove that, if $X$ is not a rosette, then $X$ arises from $\ve$. 

When $X$ is a singular subspace there is nothing to prove. So, suppose that $X$ is neither a singular subspace nor a rosette. By contradiction, suppose that $X \subset \ve^{-1}(\langle \ve(X)\rangle) =: X'$. Choose a point $x_0\in X'\setminus (X\cup X^\perp)$. As $\ve(x_0) \in \langle \ve(X)\rangle$, there exists a finite subset $A \subseteq X$ such that $\ve(x_0) \in \langle \ve(A)\rangle$.  As $X$ is not a rosette, we can always choose $A$ in such a way that $X_0 := \langle A\rangle$ contains a pair of mutually opposite lines. Hence $X_0$ is not a rosette. 

Put $X'_0 := \ve^{-1}(\langle \ve(X_0)\rangle)$. Clearly, $X'_0$ has finite rank, since $\langle \ve(X_0)\rangle$ is finite dimensional. If $X'_0$ is non-degenerate, then put $X_1 = X'_0$. Otherwise, let $R := X'_0\cap X_0'^\perp$ be the radical of $X'_0$. Then $R$ has finite rank, since $\rk(X'_0) < \infty$.  However $\cS$ is non-degenerate of infinite rank. An easy argument exploiting induction on $\rk(R)$ shows that we can find a finite subset $B$ of $\cS$ such that $\ve^{-1}(\langle\ve( X_0\cup B)\rangle)$ is non-degenerate. With $B$ chosen in this way, put $X_1 := \ve^{-1}(\langle \ve(X_0\cup B)\rangle)$. Again, $X_1$ has finite rank, since $\langle\ve(X_0\cup B)\rangle = \langle \ve(A\cup B)\rangle$ has finite dimension (less than $|A\cup B|$). Moreover $X_1$ is non-degenerate, thanks to the addition of $B$ to $A$. By construction, $X_1$ arises from $\ve$. Therefore $\ve(X_1)$ is the polar space defined on $\langle \ve(X_0\cup B)\rangle  = \langle \ve(X_1)\rangle$ by the pseudoquadratic (or alternating) form induced on $\langle \ve(X_1)\rangle$ by the pseudoquadratic (respectively, alternating) form which describes $\ve(\cS)$. This induced form is non-degenerate, since $X_1$ is non-degenerate. Moreover the embedding of $X_1$ in $\langle \ve(X_1)\rangle$, say $\ve_1$, is universal by the above and since $\ve$ is universal by assumption (which means that in the characteristic 2 case $\ve(\cS)$ cannot be described by an alternating form). Thus we can apply Theorem 1 of \cite{CGP} to $X_0$ as a subspace of $X_1$. By that theorem, $X_0$ arises from $\ve_1$, namely $X_0 = \ve_1^{-1}(\langle \ve_1(X_0)\rangle)$. However $\ve_1^{-1}(\langle \ve_1(X_0)\rangle) = \ve^{-1}(\langle \ve(X_0\rangle)$ and $x_0 \in \ve^{-1}(\langle \ve(X_0)\rangle$ by definition of $X_0$. In the end, $x_0 \in X_0$. However $X_0 \subseteq X$. Therefore $x_0 \in X$. We have reached a final contradiction.   \hfill $\Box$ 

\subsubsection{Optimally embeddable subspaces} 

Suppose that $\cS$ is embeddable and let $X$ be a subspace of $\cS$. Then $\langle \ve(X), \ve(X)^{\perp_\ve}\rangle \subseteq  \ve(X\cap X^\perp)^{\perp_\ve}$ for every embedding $\ve$ of $\cS$. Since both $\ve(X)^{\perp_\ve}$ and $\ve(X\cap X^\perp)^{\perp_\ve}$ contain the radical $R_\ve$ of $\pi_\ve$ and all quotients of $\ve$ arise by factorizing $\ve$ over subspaces of $R_\ve$, we have
\begin{equation}\label{eqR2}
\langle \ve(X), \ve(X)^{\perp_\ve}\rangle ~ = ~ \ve(X\cap X^\perp)^{\perp_\ve}
\end{equation}
if and only if the same holds for a quotient of $\ve$. Consequently, if (\ref{eqR2}) holds for $\ve$ then it also holds for all covers and all quotients of $\ve$. If $\ve(\cS) = \cS_{\pi_\ve}$ then (\ref{eqR2}) holds for every subspace $X$ such that $\dim(\ve(X)) < \infty$. (This follows from the fact that the reflexive sesquilinear form $f$ associated to $\pi_\ve$ is trace-valued and if $\pi_\ve$ is a polarity then $f$ is non-degenerate.) In particular, if $\cS$ is one of the exceptional embeddable generalized quadrangles which do not admit the absolutely universal embedding then (\ref{eqR2}) holds for every subspace $X$ and every embedding $\ve$ of $\cS$. So, if (\ref{eqR2}) holds for a subspace $X$ of $\cS$ and at last one embedding of $\cS$ then it holds for $X$ and all embeddings of $\cS$. We are now ready to state the following definition. 

\begin{defin}\label{optimal def}
\em
If $X$ satisfies property (\ref{eqR2}) for some (equivalently, every) embedding of $\cS$ then we say that $X$ is {\em optimally embeddable}. 
\end{defin}

Clearly, al singular subspaces of $\cS$ are optimally embeddable. As previously recalled, if $\cS$ admits an embedding $\ve$ such that $\ve(\cS) = \cS_{\pi_\ve}$ then every finitely generated subspace of $\cS$ is optimally embeddable. Note also that a non-degenerate subspace $X$ of $\cS$ is optimally embeddabe if and only if $\ve(X)\cup \ve(X)^{\perp_\ve}$ spans $\PG(V)$ for every embedding $\ve:\cS\rightarrow\PG(V)$ of $\cS$. 

\begin{prop}\label{optimal}
Let $Y$ and $Z$ be finite-dimensional opposite singular subspaces of $\cS$ and put $X:= \langle Y, Z\rangle$. Then $X\cap X^\perp = \emptyset$ and $\langle \ve(X), \ve(X)^{\perp_\ve}\rangle = \PG(V)$ for every embedding  $\ve:\cS\rightarrow\PG(V)$ of $\cS$.  
\end{prop} 
{\bf Proof.} Let $\ve:\cS\rightarrow\PG(V)$ be an embedding of $\cS$. The projective subspace $\langle \ve(X)\rangle$ is the union of the lines of $\PG(V)$ joining a point of $\ve(Y)$ with a point of $\ve(Z)$. Let $\ell$ be one of these lines, say $\ell = \langle \ve(y), \ve(z)\rangle$ for $y\in Y$ and $z\in Z$. Suppose that $\ell$ meets $\ve(X)^{\perp_\ve}$ non-trivially. Then at least one of $y$ and $z$ belongs to $X^\perp$. This contradicts the assumption that $Y$ and $Z$ are opposite. Therefore $\langle \ve(X)\rangle\cap \ve(X)^{\perp_\ve} = \emptyset$. Accordingly, $X\cap X^\perp = \emptyset$. 

By assumption, $\dim\langle \ve(X)\rangle < \infty$. We have $\langle {\cal X}, {\cal X}^{\perp_\ve}\rangle = ({\cal X}\cap{\cal X}^{\perp_\ve})^{\perp_\ve}$ for every finite-dimensional subspace $\cal X$ of $\PG(V)$. Therefore 
\[\langle \ve(X), \ve(X)^{\perp_\ve}\rangle ~ =  ~(\langle \ve(X)\rangle\cap\ve(X)^{\perp_\ve})^{\perp_\ve}.\]
However $\langle \ve(X)\rangle\cap\ve(X)^{\perp_\ve} = \emptyset$, as shown in the previous paragraph. Therefore  $\langle \ve(X), \ve(X)^{\perp_\ve}\rangle = \PG(V)$, as claimed.   \hfill $\Box$

\begin{lemma}\label{optimal2} 
Let $a$ and $b$ be opposite points of $\cS$. If $\ve$ is minimal then $\ve(\{a,b\})^{\perp_\ve}\cap\ve(\{a,b\})^{\perp_\ve\perp_\ve} = \emptyset$. 
\end{lemma}
{\bf Proof.} Let $\ve:\cS\rightarrow\PG(V)$ be minimal, namely $\pi_\ve$ is a polarity. Hence $\ve(\{a,b\}^{\perp_\ve}$ has codimention 2 in $\PG(V)$ and, consequently, $L := \ve(\{a,b\})^{\perp_\ve\perp_\ve}$ is a line of $\PG(V)$. By way of contradiction, suppose that $\ve(\{a,b\})^{\perp_\ve}\cap L \neq \emptyset$. Then $L\subseteq L^{\perp_\ve}$, since $L$ contains $\ve(a), \ve(b)$ and $p$, which are mutually distinct, each of these points is absolute for $\pi_\ve$ and $\pi_\ve(p)$ contains $\ve(a)$ and $\ve(b)$. Consequently $\ve(a)\perp_\ve \ve(b)$, namely $a\perp b$, while $a\not\perp b$ by assumption.  \hfill $\Box$ 

\begin{prop}\label{optimal3}
Let $a$ and $b$ be opposite points of $\cS$. Then both $\{a,b\}^\perp$ and $\{a,b\}^{\perp\perp}$ are optimally embeddable.
\end{prop}
{\bf Proof.} As $\{a,b\}^\perp\cap\{a,b\}^{\perp\perp} =  \{a,b\}^{\perp\perp\perp}\cap\{a,b\}^{\perp\perp} = \emptyset$, we need to prove that the following holds for at least one embedding $\ve:\cS\rightarrow\PG(V)$ of $\cS$.  
\begin{equation}\label{eq optimal1}
\begin{array}{rcl}
\langle \ve(\{a,b\}^\perp), \ve(\{a,b\}^\perp)^{\perp_\ve}\rangle & = & \PG(V),\\
\langle \ve(\{a,b\}^{\perp\perp}), \ve(\{a,b\}^{\perp\perp})^{\perp_\ve}\rangle & = & \PG(V).
\end{array}
\end{equation} 
We have $\langle \ve(\{a,b\}^\perp)\rangle = \{\ve(a), \ve(b)\}^{\perp_\ve}$ and $\langle \ve(\{a,b\}^{\perp\perp})\rangle = \{\ve(a), \ve(b)\}^{\perp_\ve\perp_\ve}$ by Corollary \ref{perpperp3}. Moreover $\ve(\{a,b\})^{\perp_\ve\perp_\ve\perp_\ve} = \ve(\{a,b\})^{\perp_\ve}$. Hence both equations of (\ref{eq optimal1}) are equivalent to the following:
\begin{equation}\label{eq optimal2}
\langle \ve(\{a,b\})^{\perp_\ve}, \ve(\{a,b\})^{\perp_\ve\perp_\ve}\rangle ~  = ~ \PG(V).
\end{equation} 
We can assume that $\ve$ is minimal. Hence $\ve(\{a,b\})^{\perp_\ve\perp_\ve}$ is a line, as noticed in the proof of Lemma \ref{optimal2}. Therefore
\[\langle \ve(\{a,b\})^{\perp_\ve\perp_\ve}, \ve(\{a,b\})^{\perp_\ve\perp_\ve\perp_\ve}\rangle  =  (\ve(\{a,b\})^{\perp_\ve\perp_\ve}\cap\ve(\{a,b\})^{\perp_\ve\perp_\ve\perp_\ve})^{\perp_\ve}.\]
However $\ve(\{a,b\})^{\perp_\ve\perp_\ve\perp_\ve} = \ve(\{a,b\})^{\perp_\ve}$. Hence the above equality amounts to the following:
\begin{equation}\label{eq optimal3}
\langle \ve(\{a,b\})^{\perp_\ve\perp_\ve}, \ve(\{a,b\})^{\perp_\ve}\rangle ~ = ~ (\ve(\{a,b\})^{\perp_\ve\perp_\ve}\cap\ve(\{a,b\})^{\perp_\ve})^{\perp_\ve}.
\end{equation} 
We know that $\ve(\{a,b\})^{\perp_\ve\perp_\ve}\cap\ve(\{a,b\})^{\perp_\ve} = \emptyset$ by Lemma \ref{optimal2} and, obviously, $\emptyset^{\perp_\ve} = \PG(V)$. So, (\ref{eq optimal3}) is the same as (\ref{eq optimal2}).   \hfill $\Box$ 

\section{Regularity}\label{Reg}

Throughout this section $\cS$ is a (possibly non-embeddable) non-degenerate thick-lined polar space of rank at least $2$. 

\subsection{Definitions and a proof of Theorem \ref{1.1}}\label{RDC sec} 

Lat $a$ and $b$ be two opposite points of $\cS$ and $N, N'$ generators of $\{a,b\}^\perp$. Clearly, $\{N,N'\}^\perp \supseteq \langle N\cap N', \{a,b\}^{\perp\perp}\rangle$. 

\begin{defin}\label{def1}
\em
We say that two generators $N, N'$ of $\{a,b\}^\perp$ form a {\em $\perp$-minimal pair} if $\{N,N'\}^\perp = \langle N\cap N', \{a,b\}^{\perp\perp}\rangle$. A generator $N$ of $\{a,b\}^\perp$ will be said to be {\em $\perp$-minimal} if $\{N,N'\}$ is $\perp$-minimal for every generator $N'$ of $\{a,b\}^\perp$ (in particular, $N^\perp = \langle N, \{a,b\}^{\perp\perp}\rangle$). We say that $\{a,b\}$ is {\em regular} if all generators of $\{a,b\}^\perp$ are $\perp$-minimal. If all pairs of opposite points of $\cS$ are regular, then $\cS$ is said to be {\em regular}. 
\end{defin}

\begin{lemma}\label{R0}
We have $\langle X, \{a,b\}^{\perp\perp}\rangle = \cup_{x\in \{a,b\}^{\perp\perp}}\langle X, x\rangle$ for every singular subspace $X$ of $\{a,b\}^\perp$. 
\end{lemma}
{\bf Proof.} Let $x, y$ be distinct points of $\{a,b\}^{\perp\perp}$. Then $x \not\perp y$. Consequently, no point of $\langle X, x\rangle\setminus X$ can be collinear with a point of $\langle X, y\rangle\setminus X$. The conclusion follows from this remark.  \hfill $\Box$ \\

In particular, if $N$ is a generator of $\{a,b\}^\perp$ then $\langle N, \{a,b\}^{\perp\perp}\rangle$ is the union of the generators of $\cS$ which contain $N$ and meet $\{a,b\}^{\perp\perp}$ non-trivally. Accordingly, $\{N,N\}$ is $\perp$-minimal if and only if all generators of $\cS$ containing $N$ meet $\{a,b\}^{\perp\perp}$ non-trivially. 

\begin{theo}\label{RA} 
Let $N$ be a generator of $\{a,b\}^\perp$. If the pair $\{N,N\}$ is $\perp$-minimal then $N$ is $\perp$-minimal.   

\end{theo}
{\bf Proof.} Let $\{N,N\}$ be $\perp$-minimal, namely all generators of $\cS$ containing $N$ can be obtained as $\langle N, c\rangle$ for a point $c\in \{a,b\}^{\perp\perp}$ (Lemma \ref{R0}). Given another generator $N'$ of $\{a,b\}^\perp$, let $x \in \{N, N'\}^\perp$ and suppose that $x\not\in N\cap N'$. Therefore $x\not \in N$ (otherwise $\langle x, N'\rangle$ is a singular subspace contained in $\{a,b\}^\perp$ and properly containing $N'$, a contradiction with $N'$ being a generator of $\{a,b\}^\perp$). Accordingly, $M := \langle N, x\rangle$ is a generator of $\cS$. As $\{N,N\}$ is $\perp$-minimal, $M$ contains a point $c\in \{a,b\}^\perp$. If $x = c$ then $x\in \{a,b\}^{\perp\perp}\subseteq \langle N\cap N', \{a,b\}^{\perp\perp}\rangle$. Otherwise, the line $\langle x, c\rangle$ meets $N$ in a point, say $y$. We have $\langle x, c\rangle \subseteq N'^\perp$, since both $x$ and $c$ belong to $N'^\perp$. However, $N\cap N'^\perp = N\cap N'$. Therefore $y \in N\cap N'$. Hence $x \in \langle y, c\rangle \subseteq \langle N\cap N', \{a,b\}^{\perp\perp}\rangle$. 

We have proved that $\{N, N'\}^\perp \subseteq \langle N\cap N', \{a,b\}^{\perp\perp}\rangle$. Hence $\{N, N'\}^\perp = \langle N\cap N', \{a,b\}^{\perp\perp}\rangle$ since $(N\cap N')\cup\{a,b\}^{\perp\perp} \subseteq \{N, N'\}^\perp$. So, $\{N, N'\}$ is $\perp$-minimal. As $N'$ is an arbitrary generator of $\{a,b\}^\perp$, $N$ is $\perp$-minimal.  \hfill $\Box$ \\

Theorem \ref{1.1} immediately follows from Theorem \ref{RA}. Indeed in property $(\mrR 3)$ it is assumed that every pair of generators of $\{a,b\}^\perp$ is $\perp$-minimal while $(\mrR 2)$ is equivalent to $\{N,N\}$ being $\perp$-minimal for every generator $N$ of $\{a,b\}^\perp$. Trivially, $(\mrR 3)$ implies $(\mrR 2)$. Conversely, by Theorem \ref{RA}, property $(\mrR 2)$ implies $(\mrR 3)$. 
 
\begin{cor}\label{ovvio1}
Assume that $\rk(\cS) <\infty$ and let $N$ be a generator of $\{a,b\}^\perp$. If there exists a generator $N'$ of $\{a,b\}^\perp$ such that $\{N,N'\}$ is $\perp$-minimal then $N$ is $\perp$-minimal.  
\end{cor}
{\bf Proof.} In view of Theorem \ref{RA}, we only need to prove that, if a generator $N'$ of $\{a,b\}^\perp$ exists such that $\{N,N'\}$ is $\perp$-minimal, then  $\{N,N\}$ is $\perp$-minimal. If $N' = N$ there is nothing to prove. So, assume that $N'\neq N$. Hence $N\cap N'^\perp = N\cap N'$. Let $M$ be a generator of $\cS$ containing $N$. Since $\dim(N')$ is finite and $\dim(M) = \dim(N')+1$, we have $\dim(M\cap N'^\perp) = 1+\dim(N\cap N')$ (recall that $\dim(N') = \dim(N)$ as $\rk(\cS) < \infty$). Accordingly, $M\cap N'^\perp$ is not contained in $N\cap N'$. Let $x \in (M\cap N'^\perp) \setminus (N\cap N')$. Clearly $x \in \{N, N'\}^\perp$. Therefore $x \in \langle N\cap N', \{a,b\}^\perp\rangle$ since by assumption $\{N,N'\}$ is $\perp$-minimal. By lemma \ref{R0}, the point $x$ belongs to a line $\ell$ ioining a point $c\in \{a,b\}^{\perp\perp}$ with a point $y\in N\cap N'$. We have $x \neq y$, since $x\not\in N\cap N'$. Therefore $\ell = \langle x, y\rangle$. Consequently $\ell \subseteq M$, since both $x$ and $y$ belong to $M$. However $c\in \ell$. Hence $M$ meets $\{a,b\}^{\perp\perp}$ in $c$. We have proved that all generators of $\cS$ which contain $N$ meet $\{a,b\}^\perp$ non-trivially, namely $\{N,N\}$ is $\perp$-minimal.  \hfill $\Box$

\begin{cor}\label{ovvio2}
Still assuming that $\rk(\cS) < \infty$, when $\rk(\cS) = 2$ we also assume that $\cS$ is embeddable. Suppose that there exist two opposite points $a$ and $b$ of $\cS$ such that $\{a,b\}^\perp$ admits a $\perp$-minimal pair of generators. Then $\cS$ is regular.
\end{cor}
{\bf Proof.} If $\cS$ is non-embeddable of rank $3$ then $\cS$ is regular, as proved in \cite[Proposition 5.9.4]{PVMT} (also \cite[Result 3.3]{CCGP}). Let $\cS$ be embeddable of finite rank. Then $\Aut(\cS)$ satisfies both the following transitivity properties:
\begin{itemize}
\item[$(\mathrm{T}1)$] $\Aut(\cS)$ acts transitively on the set of pairs of opposite points of $\cS$;
\item[$(\mathrm{T}2)$] for any two opposite points $x$ and $y$ of $\cS$, the stabilizer of $\{x,y\}$ in $\Aut(\cS)$ acts transitively on the set of generators of $\{a,b\}^\perp$. 
\end{itemize}
By assumption, two opposite points $a$ and $b$ exist in $\cS$ such that $\{a,b\}^\perp$ admits a $\perp$-minimal pair of generators $\{N,N'\}$. By Corollary \ref{ovvio1}, both $N$ and $N'$ are $\perp$-minimal. Hence $\{a,b\}$ is regular by $(\mathrm{T}2)$ and $\cS$ is regular by $(\mathrm{T}1)$. \hfill $\Box$ 

\begin{note}
\em
When $\rk(\cS) = 2$ and $(\mathrm{T}1)$ fails to hold (hence $\cS$ is non-embeddable) it can happen that some but not all of the pairs of opposite points of $\cS$ are regular (see \cite[Remark 6]{CCGP}). 
\end{note}

\begin{note}
\em
When $\rk(\cS)$ is infinite condition $(\mathrm{T}1)$ holds true but in general $(\mathrm{T}2)$ fails to hold. As we shall see in Section \ref{Ex more}, when $\rk(\cS) = \infty$ this case it can happen that, for any two opposite points $a$ and $b$ of $\cS$, some but not all of the generators of $\{a,b\}^\perp$ are $\perp$-minimal.  
\end{note}

\subsection{An improvement of Definition \ref{def1}}

Our definition of $\perp$-minimality as stated in Definition \ref{def1} is relative to (the hyperbolic line spanned by) a given pair $\{a,b\}$ of opposite points, but we have dropped the reference to $\{a,b\}$ in our terminology presuming that this is implicit in choosing $N$ and $N'$ among the generators of $\{a,b\}^\perp$. One might believe that, if we regard $N$ and $N'$ as sub-generators of $\cS$ without choosing in advance a pair of opposite points $\{a,b\}$ such that $\{a,b\}^\perp$ contains $N$ and $N'$, then whether $\{N,N'\}$ is or is not $\perp$-minimal depends on the choice of a pair of opposite points in  $\{N,N'\}^\perp$. However, as stated in the next theorem, this is not the case, except possibly when $\rk(\cS) = 2$ and $\cS$ is non-embeddable. So, all in all, the terminology adopted in Definition \ref{def1} is not so deceiving as it looks. 

\begin{theo}\label{def comm}
Let $N$ and $N'$ be sub-generators of $\cS$ and let $\{a,b\}$ and $\{c,d\}$ be pairs of opposite points contained in $\{N, N'\}^\perp$. Then the pair $\{N, N'\}$ is $\perp$-minimal with respect to $\{a,b\}$ if and only if it is $\perp$-minimal with respect to $\{c,d\}$, except possibly when $\cS$ is non-embeddable of rank $2$, $N = N'$ and the stabilizer in $\Aut(\cS)$ of the point $p := N = N'$ acts intransitively on the set of hyperbolic lines of $\cS$ contained in $p^\perp$. 
\end{theo} 
{\bf Proof.} Suppose firstly that $\cS$ is embeddable and let $\ve:\cS\rightarrow \PG(V)$ be an embedding of $\cS$. We can assume to have chosen $\ve$ in such a way that it is minimal, namely $\pi_\ve$ is a polarity. Then $\{a,b\}^{\perp\perp} = \ve^{-1}(\langle \ve(a), \ve(b)\rangle)$ and $\{c,d\}^{\perp\perp} = \ve^{-1}(\langle \ve(c), \ve(d)\rangle)$. Moreover, by Lemma \ref{R0} the subpspace $\langle N\cap N' , \{a,b\}^{\perp\perp}\rangle$ is the union of the subspaces $\langle N\cap N', x\rangle$ for $x\in \{a,b\}^{\perp\perp}$, A similar description holds for $\langle N\cap N', \{c,d\}^\perp\rangle$. 

Suppose that $\{N,N'\}$ is $\perp$-minimal with respect to $\{a,b\}$, namely:
\begin{equation}\label{eq def1}
\{N,N'\}^\perp ~ =~  \langle N\cap N', \{a,b\}^{\perp\perp}\rangle.
\end{equation}
By assumption, $c, d \in \{N, N'\}^\perp$ and $c\not\perp d$. Therefore, by (\ref{eq def1}) and the previous paragraph, $c\in \langle N\cap N', x\rangle\setminus N\cap N'$ and $d\in \langle N\cap N', y\rangle\setminus N\cap N'$ for suitable points $x, y \in \{a,b\}^\perp$. Without loss, we can assume that $x = a$ and $y = b$. It follows that $\langle \ve(N\cap N'), \ve(a), \ve(b)\rangle = \langle \ve(N\cap N'), \ve(c), \ve(d)\rangle$. Accordingly, $\langle \ve(c), \ve(d)\rangle$ meets every subspace $\langle \ve(N\cap N'), \ve(x)\rangle$ non-trivially, for every $x\in \{a,b\}^{\perp\perp}$. Similarly, $\langle \ve(a), \ve(b)\rangle$ meets $\langle \ve(N\cap N'), \ve(y)\rangle$ non-trivially for every $y\in \{c,d\}^{\perp\perp}$. It follows that 
\begin{equation}\label{eq def2}
\langle N\cap N', \{a,b\}^{\perp\perp}\rangle ~ = ~ \langle N\cap N', \{c,d\}^{\perp\perp}\rangle.
\end{equation}
By (\ref{eq def1}) and (\ref{eq def2}) we get $\{N,N'\}^\perp =  \langle N\cap N', \{c,d\}^{\perp\perp}\rangle$, namely $\{N,N'\}$ is also $\perp$-minimal with respect to $\{c,d\}$. 

When $\cS$ is non-embeddable of rank 3 then $\cS$ is regular  \cite[Proposition 5.9.4]{PVMT}; also \cite[Result 3.3]{CCGP}). In this case there is nothing to prove. Finally, suppose that $\rk(\cS) = 2$. Then $N'$ and $N'$ are points, say $p = N$ and $p' = N'$. If $p\neq p'$ and $\{p,p'\}$ is $\perp$-minimal with respect to $\{a,b\}$ then $\{p,p'\}^\perp = \{a,b\}^{\perp\perp}$. Consequently, since $c, d \in \{p,p'\}^\perp$ and $c \not\perp d$, we have $c, d \in \{a,b\}^{\perp\perp}$, namely $\{c,d\}^{\perp\perp} = \{a,b\}^{\perp\perp}$. Again, there is nothing to prove. Finally, let $p' = p$. If the stabilizer of $p$ in $\Aut(\cS)$ acts transitively on the set of hyperbolic lines contained in $p^\perp$ (as it is the case when $\cS$ is embeddable), then the conclusion follows.  \hfill $\Box$\\

When either $\rk(\cS) > 2$ or $\cS$ is embeddable of rank 2 Theorem \ref{def comm} allows us to replace Definition \ref{def1} with the following.    

\begin{defin}\label{def2}
\em
Let $N$ and $N'$ be sub-generators of $\cS$, possibly $N = N'$. We say that $N$ and $N'$ form a {\em $\perp$-minimal pair} if either $\{N,N'\}^\perp$ is a singular subspace or $\{N,N'\}^\perp = \langle N\cap N', \{a,b\}^{\perp\perp}\rangle$ for every (equivalently, at least one) choice of opposite points $a, b \in \{N,N'\}^\perp$. A sub-generator $N$ is said to be {\em $\perp$-minimal} if $\{N,N\}$ is $\perp$-minimal, namely either $N^\perp$ is a  generator of $\cS$ or $N^\perp = \langle N, \{a,b\}^{\perp\perp}\rangle$ for some (equivalently, every) pair $\{a,b\}$ of opposite points of $N^\perp$. A pair $\{a,b\}$ of opposite points is {\em regular} if all generators of $\{a,b\}^\perp$ are $\perp$-minimal. The polar space $\cS$ is {\em regular} if all of its sub-generators are $\perp$-minimal (equivalently, all pairs of opposite points of $\cS$ are regular). 
\end{defin}

Note that in Definition \ref{def2} when defining $\perp$-minimal sub-generators Theorem \ref{RA} is also taken into consideration. Note also that, if $N$ and $N'$ are opposite sub-generators, then not two distinct points of $\{N,N'\}^\perp$ are colliear. Hence $\{N,N'\}^\perp$ is a singular subspace if and only it is either empty or a singleton. When $\rk(\cS) < \infty$ no sub-generator is contained in one single generator and, if $N$ and $N'$ are opposite sub-generators, then $\{N,N'\}^\perp$ has the same cardinality as the set of generators containing a given sub-generator; hence $\{N,N'\}^\perp$ is never a singular subspace. 

In contrast, in general a polar space of infinite rank admits sub-generators which are contained in just one generator and pairs of opposite generators $\{N, N'\}$ such that $|\{N,N'\}^\perp| \leq 1$.    

\begin{note}
\em
When $\rk(\cS) = 2$, if $a$ and $b$ are opposite points and $c$ and $d$ are distinct (hence opposite) points of $\{a,b\}^\perp$, then $\{c,d\}^\perp = \{a,b\}^{\perp\perp}$ if and only if $\{a,b\}^\perp = \{c,d\}^{\perp\perp}$. Accordingly, $\{a,b\}$ is regular if and only if $\{c,d\}$ is regular. This is a special case of the following general fact: for any choice of singular subspaces $X_1, X_2, Y_1, Y_2$, we have $\{Y_1, Y_2\}^\perp = \{X_1, X_2\}^{\perp\perp}$ if and only if $\{X_1, X_2\}^\perp = \{Y_1, Y_2\}^{\perp\perp}$. 
\end{note}
 
\subsection{Proof of Theorem \ref{1.2 bis} and more on tight embeddings}\label{RE sec}   

Throughout this subsection the polar space $\cS$ is assumed to be embeddable. We firstly prove Theorem \ref{1.2 bis}. With the terminology introduced so far, we rephrase Theorem \ref{1.2 bis} as follows.  

\begin{theo}\label{R2}
The polar space $\cS$ admits a tight embedding if and only if 
\begin{itemize} 
\item[$(\mrR 4)$] $\cS$ admits a $\perp$-minimal pair of opposite sub-generators  $N$ and $N'$ such that $|\{N,N'\}^\perp| > 1$ and $\langle N, N'\rangle$ is optimally embeddable.  
\end{itemize}
\end{theo}
{\bf Proof.} Let $\ve:\cS\rightarrow\PG(V)$ be a minimal embedding of $\cS$. So, $\pi_\ve$ is a polarity. Accordingly, $\langle \ve(\{a, b\}^{\perp\perp})\rangle = \{\ve(a), \ve(b)\}^{\perp_\ve\perp_\ve} = \langle \ve(a), \ve(b)\rangle$ for every pair of opposite points $a, b$ of $\cS$.

Suppose that $(\mrR 4)$ holds. By assumption, $\{N,N'\}^\perp$ contains at least two points and, since $N$ and $N'$ are opposite sub-generators, all points of $\{N,N'\}^\perp$ are mutually opposite. Let $a$ and $b$ be distinct points of $\{N, N'\}^\perp$. Then $\{N, N'\}^\perp = \{a, b\}^{\perp\perp}$ since $\{N,N'\}$ is $\perp$-minimal. Moreover $\langle N, N'\rangle$ is optimally embeddable. Therefore  
\[\begin{array}{l}
\langle \ve(N), \ve(N'), \ve(a), \ve(b)\rangle  =  \langle \ve(N), \ve(N'), \ve(\{a,b\}^{\perp\perp})\rangle  = \\
= ~ \langle \ve(N), \ve(N'), \ve(\{N, N'\}^\perp)\rangle = \langle \ve(N), \ve(N'), \{\ve(N), \ve(N)\}^{\perp_\ve}\rangle  = \\
=  \PG(V). 
\end{array}\]
Hence
\begin{equation}\label{eqR3}
\langle \ve(N), \ve(N'), \ve(a), \ve(b)\rangle  ~ = ~   \PG(V). 
\end{equation}
Put $M = \langle N, a\rangle$ and $M' = \langle N', b\rangle$. Then $M$ and $M'$ are opposite generators of $\cS$. Equation (\ref{eqR3}) shows that $\langle \ve(M), \ve(M')\rangle = \PG(V)$. So, $\ve$ is tight.  

Conversely, let $\ve:\cS\rightarrow\PG(V)$ be a tight embedding of $\cS$ and let $M$ and $M'$ be generators of $\cS$ such that $\PG(V) = \langle \ve(M), \ve(M')\rangle$. Then $M\cap M' = \emptyset$, namely $M$ and $M'$ are opposite. The equality $\langle \ve(M), \ve(M')\rangle = \PG(V)$ and the fact that $M\cap M' = \emptyset$ force $\pi_\ve$ to be a polarity. Choose $a \in M$ and $b\in M'\setminus a^\perp$ and put $N = b^\perp\cap M$ and $N' = a^\perp\cap N'$. Then $N$ and $N'$ are opposite generators of $\{a,b\}^\perp$. Clearly, $\langle\langle \ve(N\cup N')\rangle, \langle \ve(a), \ve(b)\rangle\rangle = \langle \ve(M), \ve(M')\rangle$. However $\langle \ve(M), \ve(M')\rangle = \PG(V)$ by assumption and $\{N, N'\}^\perp \supseteq \{a,b\}^{\perp\perp}$. Therefore 
\begin{equation}\label{eqR3'}
\langle \ve(N), \ve(N'), \{\ve(N), \ve(N'\}^{\perp_\ve}\rangle ~ = ~  \PG(V)
\end{equation}
 and, since $\langle \ve(N), \ve(N')\rangle\cap\{\ve(N),\ve(N')\}^{\perp_\ve} = \emptyset$, necessarily 
\begin{equation}\label{eqR3''}
\{\ve(N), \ve(N')\}^{\perp_\ve} ~ = ~ \langle \ve(a), \ve(b)\rangle ~ = ~ \langle \ve(\{a,b\}^{\perp\perp})\rangle.
\end{equation}
Equality (\ref{eqR3'}) shows that $\langle N, N'\rangle$ is optimally embeddale while (\ref{eqR3''}) is equivalent to $\{N,N'\}^\perp = \{a,b\}^{\perp\perp}$, namely $\{N,N'\}$ is a $\perp$-minimal pair.  \hfill $\Box$  \\

As previously said, the condition called $(\mrR 4)$ in Theorem \ref{R2} is just a rephrasing of condition $(\mrR 4)$ of Theorem \ref{1.2 bis}. The subspaces called $X$ and $Y$ in the statement of Theorem \ref{1.2 bis} are called $N$ and $N'$ in Theorem \ref{R2}, condition (\ref{1.2 eq 1}) of $(\mrR 4)$ as stated in Theorem \ref{1.2 bis} says that $\{X,Y\}$ is $\perp$-minimal and $\{X,Y\}^\perp$ is not a singular subspace and condition (\ref{1.2 eq 2}) is equivalent to $\langle X, Y\rangle$ being optimally embeddable. Indeed, since $\{X,Y\}^\perp\not\subseteq\{X,Y\}^{\perp\perp}$, we have $\ve(X\cup Y)^{\perp_\ve} = \langle \ve((X\cup Y)^\perp)\rangle$ by Proposition \ref{perpperp1}. So, $\langle \ve(X\cup Y\cup(X\cup Y)^\perp)\rangle = \langle \ve(\langle X, Y\rangle), \ve(\langle X, Y\rangle)^{\perp_\ve}\rangle$. 

The next Corollary is just the same as Proposition \ref{1.2}. We add its proof here in order to show that Proposition \ref{1.2} indeed follows from Theorem \ref{1.2 bis}.      

\begin{cor}
Let $\cS$ be embeddable of finite rank. Then $\cS$ is regular if and only if it admits a tight embedding. 
\end{cor}
{\bf Proof.} Let $\rk(\cS) < \infty$. Then the subspace spanned by two opposite sub-generators is optimally embeddable by Proposition \ref{optimal}. Accordingly, we can drop the condition that $\langle N, N'\rangle$ is optimally embeddable from $(\mrR 4)$ of Theorem \ref{R2}. However, if we remove that condition then, in view of Corollary \ref{ovvio2}, what remains of $(\mrR 4)$ amounts to say that $\cS$ is regular.  \hfill $\Box$  

\begin{lemma}\label{RE1}
Suppose that $\rk(\cS) \geq 3$ and $\cS$ admits a unique embedding. Let $N$ and $N'$ be opposite sub-generators of $\cS$ such that the pair $\{N,N'\}$ is $\perp$-minimal and $|\{N,N'\}^\perp| > 1$. Then $\langle N, N'\rangle$ is optimally embeddable if and only if $\langle N, N'\rangle = \{a,b\}^\perp$ for every choice of distinct (hence opposite) points $a, b \in \{N,N'\}^\perp$.
\end{lemma}
{\bf Proof.} The `if' part of the statement immediately follows from the fact that $\{a,b\}^\perp$ is optimally embeddable for every choice of opposite points $a$ and $b$ (Proposition \ref{optimal3}). Note that the hypotheses that $\rk(\cS) > 2$ and $\cS$ admits a unique embedding play no role in this implication. Turning to the `only if' part, let $\ve:\cS\rightarrow\PG(V)$ be the (unique) embedding of $\cS$ and suppose that $\langle N, N'\rangle$ is optimally embeddable. Then 
\[\langle \ve(\langle N, N'\rangle), \ve(N\cup N')^{\perp_\ve}\rangle ~ = ~ \ve(\langle N, N'\rangle\cap\{N, N'\}^\perp)^{\perp_\ve}.\]
Hence, given two opposite points $a, b\in \{N, N'\}^\perp$ (which exist because $\{N, N'\}^\perp$ is non-singular by assumption), we have
\begin{equation}\label{eq RE1 1}
\langle \ve(\langle N\cup N'\rangle), \ve(\{a,b\})^{\perp_\ve\perp_\ve}\rangle ~ = ~ \PG(V) 
\end{equation}
because $\{N, N'\}^\perp = \{a,b\}^{\perp\perp}$ by assumption, $\{a,b\}^\perp\cap\{a,b\}^{\perp\perp} = \emptyset$ and $\langle \ve(\{a, b\}^{\perp\perp})\rangle = \ve(\{a,b\})^{\perp_\ve\perp_\ve}$ by Corollary \ref{perpperp3}. 

Recall that $\ve(\{a,b\})^{\perp_\ve}\cup\ve(\{a,b\})^{\perp_\ve\perp_\ve}$ spans $\PG(V)$. The embedding $\ve$ is both absolutely universal and minimal, since by assumption $\ve$ is the unique embedding of $\cS$. As $\ve$ is minimal, 
$\ve(\{a,b\})^{\perp_\ve}$ and $\ve(\{a,b\})^{\perp_\ve\perp_\ve}$ are disjoint by Lemma \ref{optimal2}. So, $\PG(V)$ is the direct sum of $\ve(\{a,b\})^{\perp_\ve}$ and $\ve(\{a,b\})^{\perp_\ve\perp_\ve}$. Moreover $\ve(\{a,b\})^{\perp_\ve} = \langle \ve(\{a,b\}^\perp)\rangle$ by Corollary \ref{perpperp3} and $\langle \ve(\langle N, N'\rangle)\rangle\subseteq \langle \ve(\{a,b\}^\perp)\rangle$, since $N\cup N' \subseteq \{a,b\}^\perp$. By comparing all this information with equation (\ref{eq RE1 1}) we see that $\langle \ve(\langle N, N'\rangle)\rangle = \langle \ve(\{a,b\}^\perp)\rangle$. 

So far we have made no use of the hypothesis that $\rk(\cS) > 2$. We shall use it now. Since $\rk(\cS) > 2$, the subspace $\langle N, N'\rangle$ is not a rosette. Therefore it arises from an embedding of $\cS$, by Proposition \ref{App1}. However $\ve$ is the unique embedding of $\cS$. Hence $\langle N, N'\rangle$ arises from $\ve$, namely $\langle N, N'\rangle = \ve^{-1}(\langle\ve(\langle N, N'\rangle)\rangle)$. Similarly, $\{a,b\}^\perp = \ve^{-1}(\langle\ve(\{a,b\}^\perp)\rangle)$. The equality $\langle \ve(\langle N, N'\rangle)\rangle = \langle \ve(\{a,b\}^\perp)\rangle$ now forces $\langle N, N'\rangle = \{a,b\}^\perp$. \hfill $\Box$ 

\begin{theo}\label{RE2}
Suppose that $\rk(\cS) \geq 3$ and $\cS$ admits a unique embedding. Then the following are equivalent:
\begin{itemize}
\item[$(1)$] the unique embedding of $\cS$ is tight;
\item[$(2)$]  $\cS$ admits admits a pair opposite sub-generators $N$ and $N'$ such that $\langle N, N'\rangle = \{a,b\}^\perp$ for two opposite points $a$ and $b$ of $\cS$;
\item[$(3)$] $\cS$ admits two opposite generators $M$ and $M'$ such that $\langle M, M'\rangle = \cS$. 
\end{itemize}  
\end{theo}
{\bf Proof.} Trivially (3) implies (1) while (1) implies (2) by Theorem \ref{R2} and Lemma \ref{RE1}. It remains to show that (2) implies (3). Given $N, N', a$ and $b$ as in (2), let $M = \langle N, a\rangle$ and $M' = \langle N',b\rangle$. Then $M$ and $M'$ are opposite generators of $\cS$. Moreover $\langle M, M'\rangle = \langle N, N', a, b\rangle = \langle \{a,b\}^\perp, a, b\rangle$, since $\langle N, N'\rangle = \{a,b\}^\perp$. However $\langle \{a,b\}^\perp, a, b\rangle = \cS$. Hence $\langle M, M'\rangle = \cS$.  \hfill $\Box$  

\begin{note}
\em 
In Lemma \ref{RE1} and Theorem \ref{RE2} we cannot drop the hypothesis that $\rk(\cS) > 2$. Indeed the sub-generators of a generalized quadrangle are its points. Hence, with $N, N', a$ and $b$ as in the hypotheses of Lemma \ref{RE1}, when $\rk(\cS) = 2$ we have $\langle N, N'\rangle = \{a,b\}^\perp$ only if $\cS$ is a grid. When $\rk(\cS) = 2$ we should replace the equality $\langle N,N'\rangle = \{a,b\}^\perp$ with $\langle \ve(\langle N, N\rangle)\rangle = \langle \ve(\{a,b\}^\perp)\rangle$, for a minimal embedding $\ve$ of $\cS$, but this condition amounts to $\langle N, N'\rangle$ being optimally embeddable.

 The hypothesis that $\cS$ admits a unique embedding is also necessary, both in Lemma \ref{RE1} and in Theorem \ref{RE2}. Indeed suppose that $\cS$ admits a tight embedding $\ve:\cS\rightarrow\PG(V)$ and let $M$ and $M'$ be opposite generators of $\cS$ such that $\langle \ve(M), \ve(M')\rangle = \PG(V)$. If $\ve$, which is minimal, is not universal, let $\vte:\cS\rightarrow \PG(\tV)$ be a proper cover of $\ve$. Then $\langle \vte(M), \vte(M')\rangle \subset \PG(\tV)$, hence $\langle M, M'\rangle \subset \cS$; similarly, $\langle\ X, Y\rangle\subset\{a,b\}^\perp$ for any two singular subspaces $X, Y\subseteq\{a,b\}^\perp$, even if $X$ and $Y$ are such that $\langle\ve(\langle X,Y\rangle)\rangle = \langle\ve(\{a,b\}^\perp)\rangle$.    
\end{note}

\begin{note}
\em 
When $\rk(\cS) < \infty$, if $\cS = \langle M, M'\rangle$ for two (necessarily opposite) generators $M$ and $M'$ then $\cS = \langle X, X'\rangle$ for any two opposite generators $X$ and $X'$. In general, as we shall see in Section \ref{Examples}, when the rank of $\cS$ is infinite $\cS$ can admit pairs of generators $\{X,X'\}$ such that $\langle X, X'\rangle = \cS$ as well as opposite generators $Y$ and $Y'$ such that $\langle Y, Y'\rangle \subset \cS$. 
\end{note}
 
\subsection{Complementary generators and regularity}\label{more}

In this subsection $\cS$ is supposed to admit the absolutely universal embedding. Hence it also admits the minimum embedding. Throughout this subsection $\ve:\cS:\rightarrow\PG(V)$ is the minimum embedding of $\cS$.  

The following definitions, which complete previously stated definitions, will help us to make our exposition more concise. Let $M$ and $M'$ be opposite generators of $\cS$. If  $\langle \ve(M), \ve(M')\rangle = \PG(V)$ then we say that $M$ and $M'$ are {\em complementary}, also that each of them is a {\em complement} of the other one. Obviously, $\cS$ admits a pair of complementary generators if and only if $\ve$ is tight.

A pair $\{N,N'\}$ of sub-generators of $\cS$ is said to be {\em optimally embeddable} if $\langle N, N'\rangle$ is optimally embeddable (Definition \ref{optimal def}); it is said to be $\perp$-{\em degenerate} if $\{N, N'\}^\perp$ is a singular subspace. Note that, according to Definition \ref{def2}, all $\perp$-degenerate pairs of opposite sub-generators are $\perp$-minimal. Note also that, if $N$ and $N'$ are opposite sub-generators, then $\{N,N'\}$ is $\perp$-degenerate precisely when $|\{N, N'\}^\perp| \leq 1$.  

\begin{lemma}\label{RR2}
No optimally embeddable pair of opposite sub-generators is $\perp$-degenerate. 
\end{lemma}
{\bf Proof.} Let $\{N, N'\}$ be an optimally embeddable pair of opposite sub-generators of $\cS$. By way of contradiction, let $|\{N,N'\}^\perp| \leq 1$. Suppose first that $\{N, N'\}^\perp = \{c\}$ for a point $c$. Let $\vte:\cS\rightarrow\PG(\tV)$ be the universal embedding of $\cS$. Since $\{N,N'\}$ is optimally embeddable, we have $\langle \vte(N), \vte(N'), \{\vte(N), \vte(N')\}^{\perp_{\vte}}\rangle = \PG(\tV)$. However, $\vte(c)$ is the unique singular point of $\{\vte(N), \vte(N')\}^{\perp_{\vte}}$. Then $\{\vte(N), \vte(N')\}^{\perp_{\vte}} \subseteq \vte(c)^{\perp_{\vte}}$. This forces $\vte(c)^{\perp_{\vte}} = PG(\tV)$, hence $c^\perp = \cS$. We have reached a contradiction. 

Therefore $\{N, N'\}^\perp = \emptyset$. Hence $\langle \vte(N), \vte(N')\rangle = \PG(\tV)$, since $\{N,N'\}$ is optimally embedded by assumption. Consequently, if $M$ is a generator of $\cS$ contaning $N$, then $\vte(M)$ meets $\vte(N')$ non-trivially. Equivalently $M$ meets $N'$ non-trivially. This contradicts the assumption that $N$ and $N'$ are opposite.  \hfill $\Box$ 

\begin{lemma}\label{RR1}
Let $\{N,N'\}$ be an optimally embeddable pair of opposite sub-generators of $\cS$ and suppose that $\{N,N'\}$ is $\perp$-minimal. Then both $N$ and $N'$ are $\perp$-minimal.
\end{lemma}
{\bf Proof.} By Lemma \ref{RR2}, $\{N,N'\}^\perp$ contains at least two distinct (hence opposite) points $a$ and $b$. So, $N$ and $N'$ are opposite generators of $\{a,b\}^\perp$ and $\{N, N'\}^\perp = \{a,b\}^{\perp\perp}$ since $\{N,N'\}$ is $\perp$-minimal. 

By way of contradiction, suppose that $N$ is not $\perp$-minimal. So, $N^\perp$ contains a point $c\not\in N\cup\{a,b\}^{\perp\perp}$. However $\{N,N'\}^\perp =\{a,b\}^{\perp\perp}$ by assumption. Moreover $\{\ve(a), \ve(b)\}^{\perp_{\ve}\perp_{\ve}} = \langle \ve(a), \ve(b)\rangle$ since $\ve$ is minimal and $\PG(V) = \langle \ve(N), \ve(N'), \ve(a), \ve(b)\rangle$, since $\langle N, N'\rangle$ is optimally embeddable. It follows that the projective plane $P = \langle \ve(a), \ve(b), \ve(c)\rangle$ meets $\langle \ve(N), \ve(N')\rangle$ in a point, say $p$. However $P \subseteq \ve(N)^{\perp_{\ve}}$. Hence $p \in X := \ve(N)^{\perp_{\ve}}\cap\langle \ve(N), \ve(N')\rangle$. However $p\not \in \ve(N)$. Therefore $X$ properly contains $\ve(N)$. Consequently, $X\cap \ve(N')\neq \emptyset$. So, $\ve(N')$ meets $X$ non-trivially. Equivalently, $N'$ meets $N^\perp$ non-trivially. This contradicts the hypothesis that $N$ and $N'$ are opposite.  \hfill $\Box$ \\

In the finite rank case the next lemma is absolutely trivial and its hypotheses are redundant. However $\rk(S) = \infty$ is allowed here. We recall that, so far, nobody has been able to prove that in every polar space of infinite rank every generator admits an opposite, although no counterexample is known which refutes this conjecture. So, the hypothesis that $M$ admits an opposite cannot be dropped from the next lemma and a proof is required for the conclusions.       

\begin{lemma}\label{RR3}
Let $M$ be a generator of $\cS$ and $p$ a point exterior to $M$. If $M$ is opposite to at least one generators of $\cS$ then there exists a generator of $\cS$ opposite to $M$ and containing $p$. If moreover $M$ admits a complement, then $M$ also admits a complement which contains $p$.  
\end{lemma}
{\bf Proof.} We firstly prove that if $M$ admits an opposite generator then a generator $M'$ also exists which is opposite to $M$ and contains $p$. 

Let $M_1$ be a generator of $\cS$ opposite to $M$. If $M_1$ contains $p$ then there is nothing to prove. Otherwise, let $M_2 := \langle p^\perp\cap M_1, p\rangle$. If $M_2\cap M = \emptyset$ then again we are done: $M' = M_2$ does the job. Assume that $M_2\cap M\neq \emptyset$. As $M_2$ contains a hyperplane $p^\perp\cap M_1$ which is disjoint from $M$, necessarily $M_2\cap M = \{a\}$ for a point $a$. Clearly, $a \neq p$. Put $N := M\cap p^\perp$. As both $p$ and $a$ belong to $N^\perp$, the line $\langle p, a\rangle$ is contained in $N^\perp$. Therefore this line meets $N$ in a point, since $N$ is a sub-generator. It follows that $a \in N$. Now, $p^\perp\not\subseteq a^\perp$ (indeed $\cS$ is non-degenerate and $p^\perp\neq a^\perp$ since $p\neq a$). Hence there exists a point $b\in p^\perp\setminus a^\perp$. Put $M_3 := \langle b , b^\perp\cap M_2\rangle$. The generator $M_3$ contains $p$. 

We claim that $M_3\cap M = \emptyset$. Indeed, suppose the contrary and let $c \in M_3\cap M$. Then $c \neq b$, as $b\not\perp a$ while $c\perp a$. The line $\langle b, c\rangle$ meets $b^\perp\cap M_2$ in a point $d$. If $c \neq d$ then $b\in \langle c, d\rangle$. However, $\langle c, d\rangle \subseteq a^\perp$, since $c \perp a$ (because $c, a \in M$) and $d\perp a$ (because $d, a \in M_2)$. It follows that $b\perp a$. This conclusion contradicts the choice of $b$. Therefore $d = c$. However $d\in b^\perp\cap M_2$, which is disjont from $M$ since $b\not\perp a$ and $a$ is the unique point of $M_2\cap M$, while $c \in M$. Again a contradiction. Hence $M_3\cap M = \emptyset$ and $M' = M_3$ is the required generator opposite to $M$ and containing $p$.  

Turning to the second claim of the lemma, suppose that $M_1$ is a complement of $M$. If $p\in M_1$ there is nothing to prove. Otherwise, let $M_2$ be constructed as above. If $M_2\cap M = \emptyset$ then $M_2$ is the required complement. Otherwise $\ve(M)\cup\ve(M_2)$ spans a hyperplane of $\PG(V)$. In this case the generator $M_3$ constructed as above is a complement of $M$.  \hfill $\Box$ 

\begin{theo}\label{RR4}
The following are equivalent:
\begin{itemize}
\item[$(1)$] $\cS$ is regular and, for every sub-generator $N$ of $\cS$, if $N^\perp$ is not a generator of $\cS$ then a sub-generator $N'$ of $\cS$ exists opposite to $N$ and such that $\{N,N'\}$ is optimally embeddable. 
\item[$(2)$] Every generator of $\cS$ admits a complement. 
\end{itemize} 
\end{theo}
{\bf Proof.} Assume (1). Let $M$ be a generator of $\cS$ and, given a point $a\in M$, let $b$ be a point opposite to $a$. Put $N = b^\perp\cap M$. Then $a, b\in N^\perp$, hence $N^\perp$ is not a generator. By (1), there exists a sub-generator $N'$ such that $\{N,N'\}$ is optimally embeddable, hence it is not $\perp$-degenerate by Lemma \ref{RR2}. However, $\cS$ is regular by (1). Hence $N^\perp = \langle N, \{a,b\}^{\perp\perp}\rangle$ and, since $\{N,N'\}$ is not $\perp$-degenerate, $N' \subset \{c,d\}^\perp$ for two opposite points $c, d\in N^\perp$. Moreover, $\{N,N'\}$ is $\perp$-minimal, namely $\{N,N'\}^\perp = \{c,d\}^{\perp\perp}$, because $\cS$ is regular. The hyperbolic line $\{c,d\}^{\perp\perp}$ meets every generator $\langle N, x\rangle$ with $x\in \{a,b\}^{\perp\perp}$. So, we can assume that $c\in \langle N,a\rangle$ and $d\in \langle N, b\rangle$. Accordingly, we can safely replace $a$ with $c$ and $b$ with $d$. In other words, we can assume that $a = c$ and $b = d$. Put $M' := \langle N', b\rangle$. Then $M$ and $M'$ are opposite. Recalling that $\{N,N'\}$ is perp-minimal and optimally embeddable, we obtain that $\ve(M)\cup\ve(M')$ spans $\PG(V)$ as in the proof of Theorem \ref{R2}. So, $M'$ is a complement of $M$.  

Conversely, assume (2). We firstly prove that every sub-generator of $\cS$ is $\perp$-minimal. Let $N$ be a sub-generator of $\cS$. If $N^\perp$ is a generator then $N$ is $\perp$-minimal (Definition \ref{def2}). Suppose that $N$ contains at least two distinct points, say $a$ and $b$. Clearly, $a\not\perp b$. Put $M := \langle N, a\rangle$. In view of (2), $M$ admits a complement. By Lemma \ref{RR3}, there exists a complement $M'$ of $M$ which contains $b$. Put $N' := a^\perp\cap M'$. As in the proof of Theorem \ref{R2}, we see that $\{N, N'\}$ is $\perp$-minimal and optimally embeddable. By Lemma \ref{RR1}, the sub-generator $N$ is $\perp$-minimal. So we have proved that every sub-generator of $\cS$ is $\perp$-minimal. Namely, $\cS$ is regular. Morover, $\{N,N'\}$ is optimally embeddable, as required in (1).   \hfill $\Box$ \\

All polar spaces constructed in Sections \ref{Ex symp}, \ref{Ex hyp} and \ref{Ex her} are regular and, in view of Proposition \ref{Q3} and its analogues for the hyperbolic and hermitian case, they satisfy condition (2) (hence (1) too) of Theorem \ref{RR4}. 

\begin{note}
\em 
When $\rk(\cS) = \infty$, the existence of a pair of complementary generators (namely the minimal embedding being tight) is not sufficient for $\cS$ to be regular and, conversely, $\cS$ being regular does not imply that all opposite generators are complementary. For instance, the quadric $\cS_{q'}$ described in Section \ref{Ex more} is not regular but it admits a pair of complementary generators, namely $[V]$ and $[V']$. Conversely, the polar spaces described in Sections \ref{Ex symp}, \ref{Ex hyp} and \ref{Ex her} are regular but each of them admits pairs of opposite generators which are not complementary. For instance, in each of these spaces, the generators called $[V]$ and $[V^*]$ are complementary. Generators also exist which are opposite to both $[V]$ and $[V^*]$. Each of them is a complement of $[V^*]$ but not of $[V]$. (See also Lemma \ref{ultima}.)
\end{note}   


\section{The three generators game}\label{3G game}  

Throughout this section $\cS$ is assumed to be embeddable. For ease of exposition, when $\rk(\cS) = 2$ we assume that $\cS$ is neither a grid of order $s > 3$ nor a generalized quadrangle as in \cite[8.6(II)(a)]{T}. So, $\cS$ admits both the absolutely universal embedding and a unique minimum embedding.   

We shall deal with triples of mutually opposite generators. Recall that a polar space $\cS$ of finite rank admits triples of opposite generators if and only if either $\cS$ is thick or $\rk(\cS)$ is even. We cannot hope for such a sharp picture for polar spaces of infinite rank, all the more that we do not even know if every generator of a polar space of infinite rank admits an opposite. However polar spaces of infinite rank exist which admit triples of mutually opposite generators. For instance, in each of the spaces described in Sections \ref{Ex symp}, \ref{Ex hyp} and \ref{Ex her}, the generators $[V]$ and $[V^*]$ are opposite and generators exist which are opposite to both $[V]$ and $[V^*]$. 
So, the theory we are going to develope in this section is not vacuous in the infinite rank case. 

\subsection{More on opposite generators}\label{more bis}  

Let $\ve:\cS\rightarrow\PG(V)$ be the minimum embedding of $\cS$. Given two opposite generators $M_1$ and $M_2$ of $\cS$, put $\cS_{M_1, M_2}:=  \ve^{-1}(\langle \ve(M_1), \ve(M_2)\rangle$. Clearly, $\cS_{M_1, M_2} \supseteq \langle M_1, M_2\rangle$. We have $\cS_{M_1, M_2} = \langle M_1, M_2\rangle$ if and only if $\ve$ is the unique embedding of $\cS$.  

\begin{lemma}\label{AddCompl0}
The following holds for $\{i,j\} = \{1,2\}$: all hyperplanes of $M_i$ are $\perp$-minimal as sub-generators of $\cS_{M_1,M_2}$ and, if $N$ is a hyperplane of $M_i$ such that $N^\perp\cap\cS_{M_1, M_2}\neq M_i$, then $N^\perp\cap M_j = \{x_j\}$ for a point $x_j\in M_j$ and $N^\perp\cap\cS_{M_1, M_2} = \langle M_i, x_j \rangle = \langle N, x_i, x_j\rangle$ for any $x_i \in M_i\setminus N$. 
\end{lemma}
{\bf Proof.} Put $\cS' := \cS_{M_1, M_2}$ for short and, denoted by $V'$ the subspace of $V$ corresponding to $\langle \ve(M_1), \ve(M_2)\rangle$, let $\ve':\cS'\rightarrow\PG(V') = \langle \ve(M_1), \ve(M_2)\rangle$ be the embedding of $\cS'$ induced by $\ve$. The embedding $\ve'$ is minimal and $\pi_\ve$ induces $\pi_{\ve'}$ on $\PG(V')$. Explicitly, $X^{\perp_{\ve'}} = X^{\perp_\ve}\cap\PG(V')$ for every subset $X$ of $\PG(V')$.

Let $N$ be a hyperplane of $M_1$ such that $M_1 \subset N^\perp\cap\cS'$. Then $\ve(N)^{\perp_{\ve'}}$ properly contains $\ve(M_1)$. Hence it meets $\ve(M_2)$ non-trivially. The intersection $\ve(N)^{\perp_{\ve'}}\cap\ve(M_2)$ cannot contain a line, otherwise $N^\perp\cap M_2$ contains a line and $\langle N, N^\perp\cap M_2\rangle$ is a singular subspace containg $N$ as a subspace of codimension at least $2$, contradcting the fact that $N$ is a sub-generator. Therefore $\ve(N)^{\perp_{\ve'}}\cap\ve(M_2) = \{\ve(x_2)\}$ for a point $x_2 \in M_2$ and $\ve(N)$ has codimension 2 in $\ve(N)^{\perp_{\ve'}}$. This implies that $\ve(N)^{\perp_{\ve'}} = \langle \ve(N), \ve(x_1), \ve(x_2)\rangle = \langle \ve(N), \{\ve(x_1), \ve(x_2)\}^{\perp_{\ve'}\perp_{\ve'}}\rangle$ for $x_1\in M_1\setminus N$. Accordingly, $N^\perp = \langle N, x_1, x_2\rangle = \langle N, (\{x_1, x_2\}^\perp\cap\cS')^\perp\cap \cS' \rangle$.  \hfill $\Box$

\begin{cor}\label{AddCompl1}
The generators $M_1$ and $M_2$ are complementary if and only if the following holds for $\{i,j\} = \{1,2\}$: 
\begin{itemize}
\item[$(\ast)$] all hyerplanes of $M_i$ are $\perp$-minimal as sub-generators of $\cS$ and, if $N$ is a hyperplane of $M_i$ such that $N^\perp\neq M_i$, then $N^\perp\cap M_j\neq \emptyset$ and $N^\perp = \langle M_i, N^\perp\cap M_j\rangle$. 
\end{itemize}
\end{cor}
{\bf Proof.} The `only if' part is Lemma \ref{AddCompl0} (see also the final part of the proof of that lemma). Turning to the `if' part, assume $(\ast)$ and let $z$ be any point of $\cS\setminus M_1$. Put $N := z^\perp\cap M_1$. Then $N^\perp \supset M_1$, since $z\not\in M_1$. Hypothesis $(\ast)$ now implies that $N^\perp\cap M = \{y\}$ for a point $y\in M_2$ and $N^\perp = \langle N, \{x,y\}^{\perp\perp}\rangle$ for $x\in M_1\setminus N$. In the first case $z \in M_2$. In the second case, $x\in \langle N, z'\rangle$ for some $z'\in \{x,y\}^{\perp\perp}$. However $\ve(\{x,y\}^{\perp\perp}) \subseteq \langle \ve(x), \ve(y)\rangle$ because $\ve$ is minimal and $\langle \ve(x), \ve(y)\rangle  \subseteq \langle \ve(M_1), \ve(M_2)\rangle$. Hence $z\in \langle \ve(M_1), \ve(M_2)\rangle$ for every point $z\in \cS$. Therefore $\langle \ve(M_1), \ve(M_2)\rangle = \PG(V)$.  \hfill $\Box$  

\begin{note}
\em
In the proof of the `if' part of Corollary \ref{AddCompl1} we have exploited only the hypothesis that $(\ast)$ holds for $(i,j) = (1,2)$. This is enough to obtain that $M_1$ and $M_2$ are complementary. Hence if $(\ast)$ holds for $(i,j) = (1,2)$ then it also holds for $(i,j) = (2,1)$.
\end{note}

\subsection{Deep and hyperbolic sub-generators} 

Let $X$ be a sub-sub-generator of $\cS$, namely a singular subspace of $\cS$ of codimension 2 in at least one (hence all) of the generators of $\cS$ which contain it. The star $\St(X)$  of $X$ is either a line (if $X$ is contained in just one generator) or a pencil of lines (if just one of the sub-generators containing $X$ is contained in at least two generators) or a non-degenerate generalized quadrangle (if none of the sub-generators containing $X$ is contained in just one generator). 

\begin{lemma}\label{easy}
If a generator $M$ of $\cS$ contains a sub-generator which is contained in just two generators, then all sub-generators contained in $M$ are contained in at most two generators. 
\end{lemma}
{\bf Proof.} By way of contradiction, suppose that $M$ contains two sub-generators $N$ and $N'$ with $N$ contained in just two generators and $N'$ in at least three generators. Put $X = N\cap N'$. Then $\St(X)$ is a non-degenerate thick-lined generalized quadrangle. However, the point of $\St(X)$ corresponding to $N$ belongs to just two lines of $\St(X)$ while the point corresponding to $N'$ belongs to at least three lines. This is impossible in any thick-lined generalized quadrangle.  \hfill $\Box$ 

\begin{defin}\label{def3}  
\em
Let $N$ be a sub-generator of $\cS$. We say that $N$ is {\em deep} if $N^\perp$ is a generator, namely $N$ is contained in just one generator; in other words, $\{N,N\}$ is $\perp$-degenerate. If $N$ is contained in exactly two generators then we say that $N$ is {\em hyperbolic}. So, by Lemma \ref{easy}, if a generator contains a hyperbolic sub-generator then it is hyperbolic (Definition \ref{1.3 def hyp}).  If all non-deep sub-generators of $\cS$ are hyperbolic (equivalently, all of its generators are hyperbolic) then we say that $\cS$ is {\em hyperbolic}. 
\end{defin}

Obviously, hyperbolic sub-generators are $\perp$-minimal. If $\cS$ admits a hyperbolic sub-generator, then at least one of the hyperbolic lines of $\cS$ has size 2. Hence all of them have size 2, since $\Aut(\cS)$ is transitive on the set of pairs of opposite points of $\cS$. So, if $\cS$ is hyperbolic then it is regular.   

If $M$ is a generator and $p$ a point exterior to $M$ then $p^\perp\cap M$ is a non-deep sub-generator. So, every generator contains non-deep sub-generators. Accordingly, a generator is non-hyperbolic if and only if all of its non-deep hyperplanes are non-hyperbolic. 

We recall that if $\rk(\cS) < \infty$ then $\cS$ admits no deep sub-generators. In this case either $\cS$ is hyperbolic or none of its sub-generators is hyperbolic (hence $\cS$ is thick). In contrast, in general, when $\rk(\cS) = \infty$ deep sub-generators exist in $\cS$. So, in general, a generator of $\cS$ contains both deep and non-deep sub-generators. Moreover, polar spaces of infinite rank also exist where some but not all of the non-deep sub-generators are hyperbolic. The quadrics discussed in Section \ref{Ex more} have this property.  

\begin{note}
\em
When $\rk(\cS) = \infty$ the stabilizer in $\Aut(\cS)$ of a generator $M$ of $\cS$ acts intransitively on the set of hyperplanes of $M$. Indeed $M$ is an infinite-dimensional projective space and the automorphism group of an infinite-dimensional projective space is never transitive on the set of hyperplanes of that space. In view of this fact, the existence of generators containing both deep and non-deep hyperplanes is not so surprising. 
\end{note} 

\subsection{An existence result for triples of opposite generators} 

The next theorem is an anlogue of Lemma \ref{RR3}, with triples of opposite generators instead of pairs. 

\begin{theo}\label{three opp}
Let $M$ and $M'$ be opposite generators of $\cS$ such that a generator of $\cS$ also exists which is opposite to both $M$ and $M'$. Suppose moreover that neither $M$ nor $M'$ are hyperbolic. Then for every point $p\not\in M\cup M'$ there exists a generator of $\cS$ opposite to both $M$ and $M'$ and containing $p$.
\end{theo}
{\bf Proof.} Let $M_1$ be a generator opposite to both $M$ and $M'$. If $p\in M_1$ there is nothig to prove. Otherwise, put $M_2 := \langle p^\perp\cap M_1, p\rangle$. If $M_2$ is opposite to both $M$ and $M'$ then we are done. Suppose firstly that $M_2$ is opposite to neither $M$ nor $M'$. Then $M\cap M_2 = \{a\}$ and $M'\cap M_2 = a'$ for distinct points $a, a' \in M_2\setminus(p^\perp\cap M_1)$. Pick a point $b \in p^\perp\setminus(a^\perp\cup a'^\perp)$ and put $M_3 := \langle b^\perp\cap M_2, b\rangle$. Then $M_3$ is a generator and, as in the proof of Lemma \ref{RR3}, we can see that $M_3$ is opposite to both $M$ and $M'$. 

Assume now that $M_2$ is opposite to just one of $M$ and $M'$, say $M_2\cap M = \emptyset$ but $M_2\cap M' := \{a\}$ for a point $a\in M_2\setminus (p^\perp\cap M_1)$. Put $X := M\cap\{p,a\}^\perp$. Then $X$ has codimension 2 in $M$, $\langle X, p, a\rangle$ is a generator and $M\cap\langle X, p, a\rangle = X$. The star $\St(X)$ of $X$ in $\cS$ is a non-degenerate polar space of rank 2 and $M$ and $\langle X, p, a\rangle$ are opposite lines in it. The star $\mathrm{St}(X)$ is not a grid, since $M$ is non-hyperbolic by assumption. Hence there exists a point $b\in X^\perp$ such that $b\perp p$ but $b\not\perp a$ and $\langle X, p, b\rangle\cap M = X$, namely $\{p,b\}^\perp\cap M = X$. 

Put $M_3 := \langle b^\perp\cap M_2, b\rangle$. Then $M_3$ is a generator, it contains $p$ and, as in the proof of Lemma \ref{RR3}, we can see that it is opposite to $M'$. By way of contradiction, suppose that $M\cap M_3 \neq\emptyset$. So $M\cap M_3 = \{c\}$ for a point $c\in M_3\setminus (b^\perp\cap M_2)$. So, $c\in \{p, b\}^\perp$. However, $\{p,b\}^\perp\cap M = X$. Therefore $c\in X$. Clearly, $c \neq b$. Hence $\langle b, c\rangle$ is a line and meets $b^\perp\cap M_2$ in a point, say $d$. We have $c\neq d$ because $d\in M_2$ and $M_2\cap M = \emptyset$. Therefore $b \in \langle c, d\rangle$. However $d \perp a$ because both $d$ and $a$ belong to $M_2$ while $c\perp a$ because $c\in X \subset M\cap a^\perp$. It follows that $\langle c, d\rangle\subseteq a^\perp$. Therefore $b \in a^\perp$. This contradicts the choice of $b$.   \hfill $\Box$ 

\begin{note}\label{three opp rem} 
\em
When $\rk(\cS) < \infty$ the hypothesis that $M$ and $M'$ are non-hyperbolic can be dropped from Theorem \ref{three opp}. Indeed in this case either $\cS$ admits no triples of opposite generators (when $\cS$ is hyperbolic of odd rank) or for any two opposite generators $M$ and $M'$ of $\cS$ and every point $p\not\in M\cup M'$ there exists a generator $M''$ opposite to both $M$ and $M'$ and containing $p$. So, either the main hypothesis of Theorem \ref{three opp} fails to hold (hence the statement of Theorem \ref{three opp} is logically valid) or the conclusion of Theorem \ref{three opp} holds true.
\end{note} 

\subsection{Partial dualities defined by pairs of opposite generators}\label{PD}  

Let $M_1$ and $M_2$ be opposite generators of $\cS$. The proof of the next lemma is easy. We leave it to the reader.
  
\begin{lemma}\label{3G0}
Let $X$ be a subspace of $M_1$ of codimension $2$ in $M_1$. Let $M_1^*(X)$ be the set of hyperplanes of $M_1$ containing $X$ (a line of the dual $M_1^*$ of $M_1$). Then one of the following occurs: 
\begin{itemize}
\item[$(1)$] We have $N^\perp\cap M_2 \neq \emptyset$ for every $N\in M_1^*(X)$, the subspace $X^\perp\cap M_2$ is a line of $M_2$ and the mapping $N\mapsto N^\perp\cap M_2$ is a bijection between $M_1^*(X)$ and $X^\perp\cap M_2$. 
\item[$(2)$] We have $N^\perp\cap M_2 = \emptyset$ for all but at most one of the hyperplanes $N\in M_1^*(N)$; moreover $X^\perp\cap M_2$ is either empty of a single point, according to whether $N^\perp\cap M_2 = \emptyset$ for either all or all but one of the hyperplanes $N\in M_1^*(X)$. 
\end{itemize} 
\end{lemma} 

Let now $M$ be a generator of $\cS$ opposite to both $M_1$ and $M_2$. We define a partial duality $\pi_{1,2}$ on $M$ as follows: for $x\in M$ let $X_1 := x^\perp\cap M_1$. If $X_1^\perp\cap M_2 = \emptyset$ then we put $\pi_{1,2}(x) = M$. Otherwise, $X_1^\perp\cap M_1 = \{x_2\}$ for a point $x_2\in M_2$; in this case we put $\pi_{1,2}(x) = x_2^\perp\cap M$. It follows from Lemma \ref{3G0} that $\pi_{1,2}$ is a partial duality of $M$. The partial duality $\pi_{2,1}$ is defined in the same way as $\pi_{1,2}$ but for permuting the roles of $M_1$ and $M_2$. 

\begin{lemma}\label{3G1}
Suppose that both $\pi_{1,2}$ and $\pi_{2,1}$ are non-degenerate. Then $\pi_{1,2}$ is a polarity if and only if $\pi_{1,2} = \pi_{2,1}$ (if and only if $\pi_{2,1}$ is a polarity). 
\end{lemma}
{\bf Proof.} Suppose that $\pi_{1,2}$ is a polarity. Then $x \in \pi_{1,2}(y)$ for every $y \in \pi_{1,2}(x)$. We have $\pi_{2,1}(x) = x_1^\perp\cap M$ where $x_1 = X_2^\perp\cap M_1$ and $X_2 = x^\perp\cap M_2$. Also, $\pi_{1,2}(y) = y_2^\perp\cap M$ with $y_2 = Y_1^\perp\cap M_2$ and $Y_1 = y^\perp\cap M_1$. However $x\in \pi_{1,2}^M(y)$. Therefore $y_2 \in X_2$. Accordingly, $x_1 \in Y_1$ and therefore $y \in \pi_{2,1}(x)$. Hence $\pi_{1,2}(x) \subseteq \pi_{2,1}(x)$. However both $\pi_{1,2}(x)$ and $\pi_{2,1}(x)$ are hyperplanes of $M$. Consequently, $\pi_{1,2}(x) = \pi_{2,1}(x)$. As this equality holds for every point $x\in M$, we have $\pi_{1,2} = \pi_{2,1}$. By reversing the previous argument we can also see that if $\pi_{1,2} = \pi_{2,1}$ then $\pi_{1,2}$ and $\pi_{2,1}$ are polarites.   \hfill $\Box$ 

\begin{lemma}\label{3G2}
For $x\in M$, we have $x\in \pi_{1,2}(x) \neq M$ if and only if $x$ belongs to a line of $\cS$ which meets both $M_1$ and $M_2$ non-trivially. Moreover, if $x\in \pi_{1,2}(x) \neq M$ then $\pi_{2,1}(x) = \pi_{1,2}(x)$. 
\end{lemma} 
{\bf Proof.} Suppose $x \in \pi_{1,2}(x) \neq M$ and let $X_1 = x^\perp\cap M_1$. Then $X_1^\perp\cap M_2 = \{x_2\}$ for a point $x_2\in M_2$ (because $\pi_{1,2}(x) \neq M$) and $x_2\perp x$ (because $x\in \pi_{1,2}(x)$). It follows that $\langle X_1, x, x_2\rangle$ is singular. However $X_1$ is a sub-generator. Therefore the line $\ell := \langle x, x_2\rangle$ meets $X_1$ in a point. So, $\ell$ meets both $M_1$ and $M_2$ non-trivially. Conversely,  if there exists a line $\ell$ of $\cS$ through $x$ which meets both $M_1$ and $M_2$ non trivially, then $x \in \pi_{1,2}(x) = \pi_{2,1}(x) = \ell^\perp\cap M$.  \hfill $\Box$\\

As in Section \ref{more bis}, let $\ve:\cS\rightarrow\PG(V)$ be the minimum embedding of $\cS$ and put $\cS_{M_1, M_2} = \ve^{-1}(\langle \ve(M_1), \ve(M_2)\rangle$. 

\begin{lemma}\label{3G3}
Suppose that $M \subseteq \cS_{M_1, M_2}$. Then $\pi_{1,2} = \pi_{2,1}$ and $\pi_{1,2}$ is a polarity. 
\end{lemma}
{\bf Proof.} Let $\ve(M) \subseteq\cS_{M_1, M_2}$. Lemma \ref{AddCompl0} implies that neither $\pi_{1,2}$ nor $\pi_{2,1}$ is degenerate. We shall prove that $\pi_{1,2} = \pi_{2,1}$. By Lemma \ref{3G1} this will be enough to obtain also that $\pi_{1,2}$ is a polarity. As in the proof of Lemma \ref{AddCompl0}, we put $\cS' := \cS_{M_1, M_2}$. Moreover, for a subset $X\subseteq \cS'$, we put $X^{\perp'} := X^\perp\cap\cS'$. 

Let $x\in M$, $X_1 = x^\perp\cap M_1$ and $x_2 \in M_2$ be such that $X_1^\perp\cap M_2 = \{x_2\}$. If $x_2\perp x$ then $x\in \pi_{1,2}(x)$. In this case the equality $\pi_{1,2}(x) = \pi_{2,1}(x)$ follows from Lemma \ref{3G2}. Suppose that $x_2\not\perp x$. Then by (the proof of) Lemma \ref{AddCompl0} we have that $X_1^{\perp'} = \langle X_1, h\rangle$ with $h = \{x, x_2\}^{\perp'\perp'}$ (a hyperbolic line of $\cS'$). Accordingly, and since $M_1\subseteq X_1^\perp$, the hyperbolic line $h$ meets $M_1$ in a point, say $x_1$. The sub-generator $X_2 := x^\perp\cap M_2$ is contained in $h^{\perp'}$. Hence 
$X_2 \perp x_1$, namely $X_2^\perp\cap M_1 = \{x_1\}$. However $\pi_{1,2}(x) = x_2^\perp\cap M = h^\perp\cap M = x_1^\perp \cap M = \pi_{2,1}(x)$. 
Hence $\pi_{1,2}(x) = \pi_{2,1}(x)$.        \hfill $\Box$ 

\begin{lemma}\label{3G4}
Suppose that $\pi_{1,2}$ is a polarity and admits at least one absolute point. Suppose morever that $\cS$ is defined over a division ring of characteristic different from $2$. Then $M \subseteq \cS_{M_1, M_2}$. 
\end{lemma}
{\bf Proof.} Let $\pi_{1,2}$ be a polarity. Then, denoted by $V_M$ the inderlying vector space of $M$ and by $\KK$ the underlying division ring of $\cS$, the polarity $\pi_{1,2}$ is defined by a reflexive sesquilinear form $f_M:V_M\times V_M\rightarrow \KK$. The form $f_M$ is non-degenerate, since polarities are non-degenerate reflexive partial dualities, by defnition. The form $f_M$ is trace-valued because $\chr(\KK) \neq 2$ by assumption \cite[\S 8.1.5]{T}. Hence either $f$ admits no non-zero isotropic vector or $V_M$ is spanned by the isotropic vectors of $f_M$. Therefore, if $\pi_{1,2}$ admits at least one absolute point, then the absolute points of $\pi_{1,2}$ generate $M$. By Lemma \ref{3G2}, all abosulte points of $\pi_{1,2}$ belong to $\cS_{M_1, M_2}$. Hence $M \subseteq \cS_{M_1,M_2}$.  \hfill $\Box$   

\begin{note}
\em
When $\cS$ is defined over a division ring of characteristic 2 it can happen that the form $f_M$ considered in the proof of Lemma \ref{3G4} is not trace-valued, namely the absolute points of $\pi_{1,2}$ do not span $M$. If this is the case then it might happen that $M\not\subseteq \cS_{M_1, M_2}$. 
\end{note} 

\subsection{Proof of Theorem \ref{3G intro}}\label{3G sec}

With $M$, $M_1$ and $M_2$ as in the previous subsection, in order to keep a record of $M$ in our notation henceforth we write $\pi_{1,2}^M$ and $\pi_{2,1}^M$ instead of $\pi_{1,2}$ and $\pi_{2,1}$. In the sequel we will often refer to the three generators property by the acronym $(3\mathrm{G})$. 

The next theorem is the same as Theorem \ref{3G intro}.

\begin{theo}\label{3G}
Suppose that $\cS$ is defined over a division ring $\KK$ of characteristic $\chr(\KK) \neq 2$. Let $M_1$ and $M_2$ be opposite generators of $\cS$ such that at least one generator of $\cS$ exists which is opposite both $M_1$ and $M_2$. Assume moreover that neither $M_1$ nor $M_2$ is hyperbolic. Then $\{M_1, M_2\}$ enjoys the three generators property  if and only if $\langle M_1, M_2\rangle = \cS$.   
\end{theo} 
{\bf Proof.} The `if' part is Lemma \ref{3G3}. Conversely, suppose that $\{M_1, M_2\}$ enjoys $(3\mathrm{G})$. Choose a line $\ell$ of $\cS$ meeting both $M_1$ and $M_2$ non-trivially and let $x$ be a point of $\ell$ different from $\ell\cap M_1$ and $\ell\cap M_2$. By Theorem \ref{three opp}, there exists a generator $M$ of $\cS$ containing $x$ and opposite to both $M_1$ and $M_2$. The point $x$ is absolute for $p_{1,2}^{M}$ by Lemma \ref{3G2}. Hence $M \subseteq \cS_{M_1, M_2}$ by Lemma \ref{3G4} and because $p_{1,2}^M$ is a polarity by $(3\mathrm{G})$ on $\{M_1,M_2\}$. Let now $y$ be any point of $x^\perp\setminus \cS_{M_1, M_2}$ and put $M' := y^\perp\cap M$. Then $M'$ is opposite to both $M_1$ and $M_2$ and $\pi_{1,2}^{M'}$ is a polarity by $(3\mathrm{G})$ on $\{M_1.M_2\}$. Moreover $x\in M'$ is absolute for $p_{1,2}^{M'}$ by Lemma \ref{3G2}. Lemma \ref{3G4} forces $M'\subseteq \cS_{M_1, M_2}$. In particular $y\in \cS_{M_1,M_2}$, a contradiction with the choice of $y$. It follows that $x^\perp\subseteq \cS_{M_1, M_2}$.

So far, we have proved that $x^\perp\subseteq \cS_{M_1,M_2}$ for every point $x\not\in M_1\cup M_2$ belonging to a line which meets both $M_1$ and $M_2$ non trivially. Let $X$ be the set of points $x$ as above and put $\overline{X} := \cup_{x\in X}x^\perp$. In order to finish the proof it is sufficient to prove that $\overline{X} \supseteq \cS\setminus(M_1\cup M_2)$. Let $y$ be a point exterior to $M_1\cup M_2$ and put $N_1 = y^\perp\cap M_1$ and $N_2 = y^\perp\cap M_2$. At most one of the points of $M_1\setminus N_1$ is collinear with $N_2$. Let $x_1\in M_1\setminus N_1$ be such that $x_1^\perp\cap M_2 \neq N_2$ and let $\ell$ be a line through $x_1$ meeting $M_2$ in a point $x_2\not\in N_2$. Then $y^\perp$ meets $\ell$ in a point $x\not\in M_1\cup M_2$. So $x\in X$ and therefore $y\in\overline{X}$.        \hfill $\Box$\\ 

As said in Conjecture \ref{3G conj}, we believe that the hypothesis that $\chr(\KK) \neq 2$ can be removed from Theorem \ref{3G}, provided that the conclusion is weakened as follows: the pair $\{M_1, M_2\}$ satisfies the three-generator property if and only if $M_1$ and $M_2$ are complementary. The next theorem offers a (faint) clue in favour of this conjecture.  

\begin{theo}\label{3Gbis}
Suppose that $\cS$ is regular and let $M_1$ and $M_2$ be opposite generators of $\cS$ such that at least one generator of $\cS$ exists which is opposite to both $M_1$ and $M_2$. Suppose moreover that neither $M_1$ nor $M_2$ is hyperbolic. Then $\{M_1, M_2\}$ satisfies the three-generator property if and only if $M_1$ and $M_2$ are complementary. 
\end{theo} 
{\bf Proof.} Lemma \ref{3G3} provides the `if' part of the theorem. Put $\cS' := \cS_{M_1, M_2}$. By definition, $M_1$ and $M_2$ are complementary if and only if $\cS' = \cS$. We shall prove that if $\cS' \subset \cS$ then $(3\mathrm{G})$ fails to hold for $\{M_1,M_2\}$, thus proving also the `only if' part. 

Suppose that $\cS' \subset \cS$. Pick a point $a\in \cS\setminus \cS'$ and let $N := a^\perp\cap M_1$. Suppose that $N^\perp\cap M_2\neq \emptyset$, say $N^\perp\cap M_2 = \{b\}$. As $\cS$ is regular, $N$ is $\perp$-minimal. Hence $N^\perp = \langle N, \{a,b\}^{\perp\perp}\rangle$ and $\{a,b\}^{\perp\perp}$ meets $M_1$ in a point $c$. Recall now that the embedding $\ve$ of $\cS$ which we use to define the subspace $\cS_{M_1, M_2} = \ve^{-1}(\langle \ve(M_1), \ve(M_2)\rangle_V)$ is the minimum one. Therefore $\langle \ve(\{a,b\}^{\perp\perp})\rangle = \langle \ve(a), \ve(b)\rangle$. It follows that the line $L := \langle \ve(a), \ve(b)\rangle$ meets $\langle \ve(M_1), \ve(M_2)\rangle$ in two distinct points, namely $\ve(b)$ and $\ve(c)$. Hence $L\subseteq\langle \ve(M_1), \ve(M_2)\rangle$, namely $\{a,b\}^{\perp\perp}\subseteq \cS' = \cS_{M_1, M_2}$. So, $a\in \cS'$. This contradicts the choice of $a$. 

Therefore $N^\perp\cap M_2 = \emptyset$. Consequently, if $M$ is a generator of $\cS$ containing $a$ and opposite both $M_1$ and $M_2$ (which exists by Theorem \ref{three opp}), then $\pi_{1,2}^M(a) = M$. Property $(3\mathrm{G})$ fails to hold for $\{M_1,M_2\}$. \hfill $\Box$ \\

Theorems \ref{3G} and \ref{3Gbis} show that, under suitable hypotheses on $\cS$, the three generators property boils down to forming a complementary pair. We shall see now that, given three mutually opposite generators, if they have not the same dimension (hence $\rk(\cS) = \infty$) then at least one of the three pairs we can form with them is not a complementary pair.  

\begin{lemma}\label{ultima} 
Let $\{M_1, M_2\}$ be a complementary pair of generators of $\cS$ with $\dim(M_1)\geq \dim(M_2)$. Then $\dim(M_3) \leq \dim(M_2)$ for every generator $M_3$ of $\cS$ opposite to both $M_1$ and $M_2$. If moreover $M_3$ is a complemet of $M_1$ then $\dim(M_3) = \dim(M_2)$. 
\end{lemma}
{\bf Proof.} Suppose that $\rk(\cS) = \infty$ (otherwise there is nothing to prove). So, the hupothesis that $M_2$ is a complement of $M_1$ combined with the hypothesis that $\dim(M_1)\geq \dim(M_2)$ imply that $\dim(M_1) = \dim(\ve) = \rk(\cS)$, with $\ve:\cS\rightarrow\PG(V)$ the minimum embedding of $\cS$. Hence $\dim(M_3) \leq \dim(M_1)$. Therefore, if $\dim(M_2) = \dim(M_1)$ then $\dim(M_3) \leq \dim(M_2)$.

Let now $\dim(M_2) < \dim(M_1)$ and, for a contradiction, suppose that $\dim(M_3) > \dim(M_2)$. Recall $\PG(V) = \langle \ve(M_1), \ve(M_2)\rangle$ since $M_1$ and $M_2$ are complementary. Let $p^{M_1}_{M_2}:\PG(V)\rightarrow\ve(M_2)$ be the projection of $\PG(V)$ onto $\ve(M_2)$ with $\ve(M_1)$ as the kernel. As $\dim(M_3) > \dim(M_2)$, $p^{M_1}_{M_2}$ induces a non-injective mapping on $\ve(M_3)$. Hence there exist points $z \in M_2$ and $x, y \in M_1$ such that $x \neq y$ and both lines $\ell_x := \langle \ve(z), \ve(x)\rangle$ and $\ell_y := \langle \ve(z), \ve(y)\rangle$ meet $\ve(M_3)$ non-trivially, say $\ell_x\cap \ve(M_3) = \{p_x\}$ and $\ell_y\cap \ve(M_3) = p_y$. Clearly $p_x \neq p_y$, since otherwise $p_x = p_y = \ve(z)$, which is impossible since $M_3\cap M_2 = \emptyset$. Therefore $\ell := \langle p_x, p_y\rangle$ is a line and $\ell$ meets $\langle \ve(x), \ve(y)\rangle$ in a point, say $p$. This point belongs to both $\ve(M_3)$ and $\ve(M_1)$. This conclusion contradict the hypothesis that $M_3$ is opposite to $M_1$. Therefore $\dim(M_3) \leq \dim(M_2)$.

When $M_3$ is also a complement of $M_1$ we can permute $M_2$ and $M_3$ in the previous paragraph, thus obtaining that $\dim(M_2) \leq \dim(M_3)$. In this case $\dim(M_3) = \dim(M_2)$.   \hfill $\Box$

\begin{cor}
Let $M_1, M_2$ and $M_3$ be three mutually opposite generators of $\cS$ and suppose that they have not the same dimension. Then at least one of $\{M_1, M_2\}, \{M_1, M_3\}$ or $\{M_2, M_3\}$ is not a complementary pair.   
\end{cor} 
{\bf Proof.} Since $M_1, M_2, M_3$ have not the same dimension, necessarily $\rk(\cS) = \infty$. To fix ideas, suppose that both $M_2$ and $M_3$ are complements of $M_1$. Suppose firstly that $\dim(M_1) \geq \dim(M_2)$. Then $\dim(M_2) = \dim(M_3)$ by Lemma \ref{ultima}. By assumption, $M_1, M_2$ and $M_3$ have not the same dimension. Hence $\dim(M_1) > \dim(M_2) = \dim(M_3)$. It follows that $\dim(M_1) = \dim(\ve) = \rk(\cS)$, where $\ve:\cS\rightarrow\PG(V)$ is the minimum embedding of $\cS$. Hence $\dim\langle \ve(M_2), \ve(M_3)\rangle = \dim(M_2) = \dim(M_3) < \dim(M_1) = \dim(\ve)$. So, $\langle \ve(M_1), \ve(M_2)\rangle \subset \PG(V)$. The generators $M_3$ is not a compement of $M_2$. 

On the other hand, let $\dim(M_1)$ be less than $\dim(M_2)$ or $\dim(M_3)$, say $\dim(M_1) < \dim(M_2)$. Then $\dim(\ve) = \dim(M_2) = \mathrm{max}(\dim(M_1), \dim(M_3))$, since both pairs $\{M_1, M_2\}$ and $\{M_1, M_3\}$ are complementary. As $\dim(M_1) < \dim(M_2)$, the equality $\dim(M_2) = \mathrm{max}(\dim(M_1), \dim(M_3))$ implies that $\dim(M_3) = \dim(M_2)$. So, $\dim(M_1) < \dim(M_2) = \dim(M_3)$. If $M_3$ is a complement of $M_2$  then, by permuting $M_1$ and $M_3$ in the statement of Lemma \ref{ultima}, we obtain that $\dim(M_1) = \dim(M_2)$, while $\dim(M_1) < \dim(M_2)$ by assumption. Therefore $\{M_2, M_3\}$ is not a complementary pair. (Indeed now we would obtain that $M_2\cap M_3 \neq \emptyset$, thus contradicting the hypothesis that $M_1, M_2$ and $M_3$ are mutually opposite.)  \hfill $\Box$ 


\section{A few examples of infinite rank}\label{Examples} 

Throughout this section $\FF$ is a given field, $V$ is a infinite dimensional $\FF$-vector space, $V^*$ is the dual of $V$ and $\ovV:= V\oplus V^*$. Recall that $\dim(V) < \dim(V^*)$. So $\dim(\ovV) = \dim(V^*)$.

We will adopt the following notation. We denote by $[v]$ the point of $\PG(\ovV)$ represented by a non-zero vector $v$ of $\ovV$ and by $[W]$ the subspace of $\PG(\ovV)$ corresponding to a subspace $W$ of $\ovV$. 

The kernel of a linear functional $\xi\in V^*$ is denoted by $\ker(\xi)$ and, for a subspace $\Xi$ of $V^*$, we put $\ker(\Xi) = \cap_{\xi\in\Xi}\ker(\xi)$. Given $x \in V$, we define $\ker^*(x) := \{\xi\in V^*~|~ \xi(x) = 0\}$ (a subspace of $V^*$) and, for a subspace $X$ of $V$, we put $\ker^*(X) := \cap_{x\in X}\ker^*(x)$. 

Finally, $p_V$ and $p_{V^*}$ are the projections of $\ovV$ onto $V$ and $V^*$ respectively with respectively $V^*$ and $V$ as the kernels and, for a subspace $W$ of $\ovV$, $p_V^W$ and $p_{V^*}^W$ are their restrictions to $W$.   

All polar spaces to be considered in this section are full subgeometries of $\PG(\ovV)$. In order to avoid confusions between generation in a polar space $\cS\subseteq \PG(\ovV)$ and spans in $\PG(\ovV)$, henceforth we renounce to use the symbol $\langle .\rangle$ when referring to spans in $\PG(\ovV)$, keeping it for generation in $\cS$ and using phrases as ``the span of ... in $\PG(\ovV)$" or switching to spans in $\ovV$ when dealing with spans in $\PG(\ovV)$. However we keep the symbol $\langle.\rangle$ for spans in $\ovV$. No ambiguity arises from doing so, since the vectors of $\ovV$ are not points of $\PG(\ovV)$, let alone points of $\cS$. Thus, if $X\subseteq \ovV$ then $\langle X\rangle$ is the span of $X$ on $\ovV$, $[\langle X\rangle]$ is the span of $[X] := \{[x]\}_{x\in X\setminus\{0\}}$ in $\PG(\ovV)$ and, if $[X]\subseteq \cS$, then $\langle [X]\rangle$ is the subspace of $\cS$ generated by $[X]$. 

Turning to perps, let $f$ be the bilinear or hermitian form which defines the quasi-polarity of $\PG(\ovV)$ associated to $\cS$ (or the quasi-polarity of $\PG(\ovV')$ associated to $\cS$ if $\cS$ lives in $\PG(\ovV')$ for a subspace $\ovV'$ of $\ovV$). We use the symbol $\perp_f$ for orthogonality with respect to $f$ both in $\PG(\ovV)$ and $\ovV$ (respectively, in $\PG(\ovV')$ and $\ovV'$), keeping the symbol $\perp$ with no subscript for collinearity in $\cS$. 

\subsection{A symplectic polar space}\label{Ex symp}

As in \cite{Pinfty}, we define a (non-degenerate) alternating form $f:\ovV\times\ovV\rightarrow \FF$ as follows: 
\begin{equation}\label{eq f}
f(a\oplus \alpha, b\oplus\beta) ~=~ \alpha(a)-\beta(b), \hspace{5 mm} \forall a,b \in V, ~\alpha, \beta \in V^*.
\end{equation}
Let $\cS_f$ be the symplectic polar space defined by $f$ in $\PG(\ovV)$. Of course, the inclusion maping of $\cS_f$ into $\PG(\ovV)$ is the minimum embedding of $\cS_f$ (the unique embedding of $\cS_f$ when $\chr(\FF) \neq 2$.). This embedding is tight. Indeed the subspaces $V = V\oplus\{0\}$ and $V^* = \{0\}\oplus V^*$ of $\ovV$ yield opposite generators $[V]$ and $[V^*]$ of $\cS_f$ and $V+V^* = \ovV$ (but if $\chr(\FF) = 2$ then $\langle [V], [V^*]\rangle \subset \cS_f$). Clearly, $\rk(\cS_f) = \dim(V^*)$.  

As proved in \cite{CCGP}, all symplectic polar spaces are regular. Hence $\cS_f$ is regular. So, if $N$ is a non-deep sub-generator of $\cS_f$ then $N^\perp$ contains a hyperbolic line $h$ and $\{\langle N, x\rangle\}_{x\in h}$ is the set of generators containing $N$. Deep sub-generators also exist in $\cS_f$. For instance, $V^*$ admits infinitely many hyperplanes $H$ such that $H \neq \ker^*(x)$ for any $x\in V$. If $N = [H]$ for such a hyperplane $H$ then $N^\perp = [V^*]$.    

\subsubsection{Singular subspaces and generators of $\cS_f$}\label{Ex symp gen} 

Let $W$ be a subspace of $\ovV$. Put $K_1 := \ker(p_{V^*}^W) = W\cap V$ and $K_2 := \ker(p_V^W) = W\cap V^*$ and let $C$ be a complement of $K_1\oplus K_2$ in $W$. Then $p_V^C$ and $p_{V^*}^C$ are injective. Accordingly, if $C_1 := p_V^C(C)$ and $C_2 := p_{V^*}^C(C)$, then $C_1 \cong C \cong C_2$ and the mappings $p_{1,2} := p_{V^*}^C\circ (p_V^C)^{-1}$ and $p_{2,1} = p_{V^*}^C\circ(p_V^C)^{-1}$ are isomorphisms from $C_1$ to $C_2$ and from $C_2$ to $C_1$ respectively. Also, $p_{2,1} = p_{1,2}^{-1}$. So, $C = \{x \oplus p_{1,2}(x)\}_{x\in C_1} = \{p_{2,1}(\xi)\oplus \xi\}_{\xi\in C_2}$. 

\begin{lemma}\label{Q0} 	
With $C_1, C_2, K_1, K_2$ as above, $C_i\cap K_i = \{0\}$ for $i = 1, 2$.
\end{lemma}
{\bf Proof.} If $x\in C_1\cap K_1$ then $W$ contains $(x\oplus p_{1,2}(x))-x = p_{1,2}(x)$. Therefore $p_{1,2}(x)\in W\cap V^* =K_2$, namely $x\oplus p_{1,2}(x) \in K_1\oplus K_2$. Hence $x\oplus p_{1,2}(x) = 0$, namely $x = 0$.  \hfill $\Box$\\

Clearly, $[W]$ is a singular subspace of $\cS_f$ if and only if 
\begin{equation}\label{eq Q0}
\left.\begin{array}{lll}
K_1 \subseteq \ker(C_2+K_2), & K_2\subseteq \ker^*(C_1+K_1) & \mbox{and}\\ 
p_{1,2}(x)(y) = p_{1,2}(y)(x), & \forall x, y\in C_1. & 
\end{array}\right\}
\end{equation} 
where for subspaces $X\subseteq V$ and $\Xi\subseteq V^*$ we put 
\[\ker^*(X) := \{\xi\in V^*~|~ \ker(\xi) \supseteq X\} \mbox{ and } \ker(\Xi) :=  \cap_{\xi\in\Xi}\ker(\xi).\] 

\begin{lemma}\label{Q1}
Assume that $[W]$ is a singular subspace of $\cS_f$. Then $[W]$ is a generator of $\cS_f$ if and only if
\begin{equation}\label{eq Q1}
\ker^*(C_1+K_1) = K_2 ~\mbox{ and } ~\ker(C_2+K_2) = K_1.
\end{equation}
\end{lemma} 
{\bf Proof.} Suppose that $[W]$ is a generator of $\cS_f$. Clearly $W^{\perp_f}$ contains  both $\ker(C_2+K_2)$ and $\ker^*(C_1+K_1)$. However $W$ is maximal among the totally $f$-isotropic subspaces of $\ovV$. Hence both $\ker(C_2+K_2)$ and $\ker^*(C_1+K_1)$ are contained in $W$, namely $\ker(C_2+K_2) = K_1$ and $\ker^*(C_1+K_1) = K_2$, as claimed in (\ref{eq Q1}). 

Conversely, assume (\ref{eq Q1}). Note that, if $X$ is a subspace of $V$, then for every vector $x\in V\setminus X$ there exists a hyperplane $H$ of $V$ which contains $X$ but does not contain $x$. Consequently, 
\begin{equation}\label{eq kerker}
\ker(\ker^*(X)) ~ = ~ X.
\end{equation}
By (\ref{eq kerker}) and the first equality of (\ref{eq Q1}) we obtain 
\begin{equation}\label{eq Q2}
C_1+K_1 ~ =~  \ker(K_2).
\end{equation}
Let now $a\oplus \alpha \in W^{\perp_f}$. Then $0 = f(a\oplus\alpha, \kappa) = \kappa(a)$ for every $\kappa\in K_2$. Therefore $a\in \ker(K_2)$. Hence $a = c+k$ for some $c\in C_1$ and some $k\in K_1$ by (\ref{eq Q2}). So, we can replace $a\oplus \alpha$ with $c\oplus\alpha = (a\oplus\alpha)-(k\oplus 0)$. In other words, we can assume to have chosen $a\oplus\alpha$ with $a\in C_1$. However $W$ contains $a\oplus p_{1,2}(a)$. Hence we can replace $a\oplus\alpha$ with $a\oplus\alpha - a\oplus p_{1,2}(a) = 0\oplus(\alpha-p_{1,2}(a)) = \alpha-p_{1,2}(a)$. So, we can assume to have chosen $a\oplus\alpha$ with $a = 0$. Thus, we are reduced to the case of $\alpha\in V^*$ such that $\alpha \perp_f W$. This condition forces $\alpha\in \ker^*(C_1+K_1)$. Hence $\alpha \in K_2$ by the first equality of (\ref{eq Q1}). We have proved that $W^{\perp_f} = W$. Thus $[W$] is a generator of $\cS_f$. \hfill $\Box$ 

\begin{prop}\label{Q3}
Every generator of $\cS_f$ admits a complement.  
\end{prop}
{\bf Proof.}  Let $[W]$ be a generator of $\cS_f$. With $C_1, C_2, K_1, K_2$ defined as in the first paragraph of this subsection, choose a complement $H_1$ of $C_1+K_1$ in $V$ and put $H_2 := \ker^*(H_1)$. Let $W' := H_1\oplus H_2$. The two relations of (\ref{eq Q1}), referred to $W' = H_1\oplus H_2$, amount to the equality  $\ker^*(H_1) = H_2$, which holds by definition of $H_2$, and the equality $\ker(H_2) = H_1$, namely $\ker(\ker^*(H_1)) = H_1$, which is a special case of (\ref{eq kerker}). So, $W'$ satisfies both equations of (\ref{eq Q1}). Hence $[W']$ is a generator.  

Let $x\oplus\xi\in W\cap W'$. Then $x\in H_1$, $\xi\in H_2$ and $x\oplus\xi = (k+c)\oplus(\kappa + p_{1,2}(c))$ for suitable $k\in K_1$, $\kappa\in K_2$ and $c\in C_1$. However $H_1\cap (K_1+C_1) = 0$. Therefore $x = 0$. Accordingly, $k+c = 0$, namely $c = -k$. Lemma \ref{Q0} now forces $c = k = 0$. Hence $p_{1,2}(c) = 0$. In the end, $\xi = \kappa$. Recall that $\xi \in H_2 = \ker^*(H_1)$ while $\kappa \in K_2 = \ker^*(K_1+C_1)$. Moreover \[\ker^*(H_1)\cap \ker^*(K_1+C_1) ~ = ~ \ker^*(H_1+(K_1+C_1)).\] 
However $H_1+K_1+C_1 = V$ by our choice of $H_1$. Hence $H_2\cap K_2 = \ker^*(V) = 0$. It follows that $\xi = \kappa = 0$. Therefore $W\cap W' = 0$. 

We shall now prove that $W+W' = \ovV$. The subspace $W+W'$ contains $K_2+H_2$. Let $\xi\in V^*$. If $\ker(\xi)$ contains $H_1$ then $\xi$ belongs to $H_2$. Suppose that $\ker(\xi)\cap H_1$ is a hyperplane of $H_1$. As $H_1$ is a complement of $K_1+H_1$ and $K_2 = \ker^*(K_1+H_1)$, there exists $\kappa \in K_2$ such that $\kappa$ and $\xi$ induce the same linear functional on $H_1$. Put $\chi := \xi-\kappa$. Then $\chi \in \ker^*(H_1) = H_2$. So, $\xi = \kappa + \chi \in K_2+H_2$. 

We have proved that $W+W'$ contains $V^*$. Consequently $(W+W')\cap V$ contains $(K_1+C_1)+H_1 = V$. So, $W+W'$ contains both $V$ and $V^*$. Hence $W+W'  = \ovV$.  \hfill $\Box$  \\

The space $\cS_f$ also admits opposite generators which are not complementary, as in the following example.  

\begin{ex}\label{esempio1}
\em
Given a basis $(e_i)_{i\in \mathfrak{I}}$ of $V$, for every $i\in \mathfrak{I}$ define a linear functional $\eta_i\in V^*$ by the clause that $\eta_i(e_j) = \delta_{i,j}$ (Kronecker symbol) for any $j\in \mathfrak{I}$. The linear functionals $\eta_i$ defined in this way span a proper subspace $V'$ of $V^*$ isomorphic to $V$. Note that $\ker(V') = 0$. Let $W$ be the subspace of $\ovV$ spanned by the vectors $e_i\oplus\eta_i$ for $i\in \mathfrak{I}$. It is easy to see that $[W]$ is a generator of $\cS_f$. Clearly, $V\cap W = 0$. However $\dim(V+W) = \dim(V) < \dim(V^*) = \dim(\ovV)$. Hence $V+W \subset \ovV$. So, $[W]$ is opposite to $[V]$ but it is not a complement of $[V]$. 
\end{ex}      

\subsubsection{The characteristic two case}\label{symp car 2}

The natural embeding of $\cS_f$ in $\PG(\ovV)$ is universal if and only if $\chr(\FF) \neq 2$. Accordingly, when $\chr(\FF) \neq 2$ the union of any two complementary generators generates $\cS_f$. In contrast, when $\chr(\FF) = 2$ no two generators of $\cS_f$ generate $\cS_f$, even when they are complementary. For instance, when $\chr(\FF) = 2$ the subspace $\cQ := \langle [V],[V^*]\rangle$ of $\cS_f$ is the set of points of $\PG(\ovV)$ represented by the non-zero vectors $a\oplus\alpha$ of $\ovV$ such that $\alpha(a) = 0$. The set $\cQ$ is a proper subspace of $\cS_f$. In fact $\cQ$ is an infinite dimensional analogue of a hyperbolic quadric, defined by the quadratic form which maps every vector $a\oplus\alpha \in \ovV$ onto the scalar $\alpha(a) \in \FF$. Referring to the next subsection (Section \ref{Ex hyp}) for a discussion of these quadrics, we devote the next paragraph to a description the universal embedding of $\cS_f$ when $\chr(\FF) = 2$. 

Assuming that $\chr(\FF) = 2$, let $\FF^2 := \{t^2\}_{t\in \FF}$ be the square subfield of $\FF$ and let $\FF^{1/2}$ be the extension of $\FF$ obtained by adding a square root $t^{1/2}$ of $t$ for every $t\in \FF\setminus \FF^2$. Of course, $\FF^2\subseteq \FF\subseteq \FF^{1/2}$, with $\FF^2 = \FF = \FF^{1/2}$ if and only if $\FF$ is perfect. As $\FF^{1/2}$ contains $\FF$, the field $\FF^{1/2}$ is naturally equipped with an $\FF$-vector space structure. Put $\tV := V\oplus \ovV \oplus \FF^{1/2}$ and define a quadratic form $q:\tV\rightarrow \FF$ as follows: 
\[q(x\oplus \xi\oplus t) ~= ~ \xi(x)+t^2, \hspace{5 mm} \forall x\in V, \xi\in V^*, t\in \FF^{1/2}.\]
The polar space $\cS_q$ defined by $q$ in $\PG(\tV)$ is isomorphic to $\cS_f$, the isomorphism being induced by the linear mapping from $\ovV$ to $\tV$ which maps every vector $x\oplus \xi\in \ovV$ onto the $\phi$-singular vector $x\oplus\xi\oplus\xi(x)^{1/2}\in \tV$. This isomorphism yields the universal embedding $\vte:\cS_f\rightarrow\PG(\tV)$ of $\cS_f$. The subspace $\cQ= \langle [V], [V^*]\rangle$ of $\cS_f$ is mapped by $\vte$ onto $\vte(\cS_q)\cap[\ovV\oplus\{0\}]$. 

\subsection{Hyperbolic quadrics of infinite rank}\label{Ex hyp} 

Without assuming any hypothesis on $\FF$, let $q:\ovV\rightarrow \FF$ be the non-degenerate quadratic form defined on $\ovV$ as follows: \begin{equation}\label{eq q}
q(a\oplus\alpha) ~ = ~ \alpha(a), \hspace{5 mm} \forall a\in V, ~ \alpha\in V^*.
\end{equation}
The bilinearized of $q$ is the symmetric bilinear form $f_q$ defined as follows: 
\begin{equation}\label{eq fq}
f_q(a\oplus \alpha, b\oplus \beta) ~ = ~ \alpha(b)+\beta(a), \hspace{5 mm}  \forall a,b \in V, ~\alpha, \beta \in V^*.
\end{equation} 
Let $\cS_q$ be the polar space defined by $q$. Note that $f_q$ non-degenerate, no matter what $\chr(\FF)$ is.  Hence the inclusion mapping of $\cS_q$ in $\PG(\ovV)$ is the unique embedding of $\cS_q$. As in Section \ref{Ex symp}, the subspaces $[V]$ and $[V^*]$ are opposite generators of $\cS_q$ and span $\PG(\ovV)$. Hence the embedding of $\cS_q$ in $\PG(\ovV)$ (which is the unique embedding of $\cS_q$) is tight. 

\begin{note}
\em
When $\chr(\FF) = 2$ the form $f_q$ is alternating and, with $f = f_q$, the polar space $\cS_q$ is the same as the subspace of $\cS_f$ called $\cQ$ in Section \ref{symp car 2}. As noticed at the end of Section \ref{symp car 2}, the space $\cQ$ admits an embedding in $\PG(\ovV\oplus\{0\})$ which however, since $\cS_q$ admits a unique embedding, is isomorphic to the natural embedding of $\cQ = \cS_q$ in $\PG(\ovV)$.
\end{note}      

\begin{lemma}\label{hyp lines0}
For every line $\ell$ of $\PG(\ovV)$, if $|\ell\cap\cS_q| > 2$ then $\ell\subseteq \cS_q$.
\end{lemma}
{\bf Proof.} Let $[a\oplus\alpha], [b\oplus\beta]$ and $[c\oplus\gamma]$ be three distinct points of $\cS_q$ on the same line $\ell$ of $\PG(\ovV)$. We can assume that $c\oplus\gamma = a\oplus\alpha + b\oplus\beta$. By assumption $\alpha(a) = \beta(b) = \gamma(c) = 0$. Hence $\alpha(b) + \beta(a) = 0$, namely $f_q(a\oplus\alpha, b\oplus\beta) = 0$. Therefore $[a\oplus\alpha]\perp_{f_q}[b\oplus\beta]$. Hence $\ell$ is totally $q$-singular. \hfill $\Box$  

\begin{prop}\label{hyp lines1}
The hyperbolic lines of $\cS_q$ have size 2. 
\end{prop}
{\bf Proof.} As $f_q$ is non-degenerate, every hyperbolic line of $\cS_q$ is contained in a projective line of $\PG(\ovV)$. Lemma \ref{hyp lines0} yields the conclusion.   \hfill $\Box$\\

All we have said in Section \ref{Ex symp gen} on singular subspaces of $\cS_f$ is word by word valid for singular subspaces of $\cS_q$, but for replacing the last condition of (\ref{eq Q0}) wth the following: $p_{1,2}(x) = 0$ for every $x\in C_1$. In particular, the generators of $\cS_q$ are still characterized by the two conditions of (\ref{eq Q1}). We have $[W]^{\perp_{f_q}} = [W]$ for every generator $[W]$ of $\cS_q$ and the analogue of Proposition \ref{Q3} holds: every generator of $\cS_q$ admits a complement. Note also that, since the inclusion mapping of $\cS_q$ in $\PG(\ovV)$ is the unique embedding of $\cS_q$, two opposite generators of $\cS_q$ are complementary only if their union generates $\cS_q$.  

\begin{lemma}\label{hyp lines2}
If a plane $P$ of $\PG(\ovV)$ contains three distinct lines of $\cS_q$, then $P\subseteq\cS_q$. 
\end{lemma}
{\bf Proof.} Let $P$ be a plane of $\PG(\ovV)$ containing three distinct lines of $\cS_q$. Every line of $P$ meets each of those three lines non trivially and, if meets no two of them in the same point, then it is fully contained in $\cS_q$ by Lemma \ref{hyp lines0}. Consequently, $P\subseteq \cS_q$. \hfill $\Box$

\begin{prop}\label{Q2}
All non-deep sub-generators of $\cS_q$ are hyperbolic. 
\end{prop}
{\bf Proof.} Let $N$ be a sub-generator of $\cS_q$ and suppose that $N^\perp$ contains two opposite points of $\cS_q$. Since $\Aut(\cS_q)$ acts transitively on the pairs of opposite points of $\cS_q$, we can safely assume that $N^\perp$ contains $[a]$ and $[\alpha]$ for a vector $a\in V$ and a linear functional $\alpha\in V^*$ such that $\alpha(a) = 1$. We shall prove that $\langle N, [a]\rangle$ and $\langle N, [\alpha]\rangle$ are the unique generators of $\cS_q$ which contain $N$. Put $V_0 := \ker(\alpha)$. Clearly, $\ker^*(a) \cong V_0^\ast$. Freely regarding $V^*_0$ as the same as $\ker^*(a)$, put $\ovV_0 := V_0\oplus V_0^* = \ker(\alpha)\oplus\ker^*(a)$. So, $N = [W]$ for a maximal totally $q$-singular subspace $W$ of $\ovV_0$. 

For a contradiction, let $\langle N, [b\oplus\beta]\rangle$ be a generator of $\cS_q$ containing $N$ and different from both $\langle N, [a]\rangle$ and $\langle N,[\alpha]\rangle$. The points $[a]$, $[\alpha]$ and $[b\oplus\beta]$ are mutually non collinear in $\cS_q$. Hence they span a plane $P$ of $\PG(\ovV)$ contained in $N^{\perp_{f_q}}$ and $P\cap N = \emptyset$ by Lemma \ref{hyp lines2}. However, $\{[a], [\alpha]\}^{\perp_{f_q}} = [\ovV_0]$ has codimension 2 in $\PG(\ovV)$ and the line of $\PG(\ovV)$ spanned by $[a]$ and $[\alpha]$ is a complement of $\{[a], [\alpha]\}^{\perp_{f_q}}$. Hence $P$ meets $\{[a], [\alpha]\}^{\perp_{f_q}}$ non trivially. Let $[c\oplus\gamma]$ be a point of $P\cap \{[a], [\alpha]\}^{\perp_{f_q}}$. Then $c\oplus\gamma \in W^{\perp_{f_q}}\cap\ovV_0$. As in the final part of the proof of Lemma \ref{Q1}, but with $V, V^*$ and $\ovV$ replaced by $V_0, V_0^*$ and $\ovV_0$ respecttively, we obtain that $c\oplus\gamma \in W$. This contradicts the fact that, as previously noticed, $P\cap N = \emptyset$.    \hfill $\Box$\\

Deep sub-generators also exists in $\cS_q$. For instance, as in $\cS_f$, every hyperplane $H$ of $V^*$ such that $\ker(H) = \{0\}$ yields a deep sub-generator $[H]$ of $\cS_q$. 

\begin{cor}\label{Q2bis}
The polar space $\cS_q$ is regular.
\end{cor}
{\bf Proof.} By Proposition \ref{hyp lines1}, the space $\cS_q$ is hyperbolic (Definition \ref{def3}).  All hyperbolic polar spaces are regular. \hfill $\Box$ \\

As noticed at the end of Section \ref{Ex symp gen}, opposite generators $[W]$ and $[W']$ of $\cS_f$ exist which are not complementry. The same occurs in $\cS_q$, as we can see by a slight modification of Example \ref{esempio1}. 

\begin{ex}\label{esempio2}
\em
Given a basis $(e_i)_{I\in \mathfrak{I}}$ of $V$, choose the linear functionals $\eta_i$ in such a way that $\eta_i(e_i) = 0$ for every $i$ and $\eta_i(e_j)+\eta_j(e_i) = 0$ for any choice of $i, j \in \mathfrak{I}$. To this goal we can use the following trick. Let $P$ be a partition of $\mathfrak{I}$ in subsets of size 2 and define two sets $P^+$ and $P^-$ of ordered pairs as follows: for every pair $\{i,j\}\in P$, choose one of the two ordered pairs $(i,j)$ and $(j,i)$ and put it in $P^+$, putting the other one in $P^-$. Then define $\eta_j(e_i) = 0$ if $\{i,j\}\not\in P$ (in particular, if $i = j$), $\eta_j(e_i) = 1$ if $(i,j)\in P^+$ and $\eta_j(e_i) = -1$ if $(i,j)\in P^-$. 

As in Example \ref{esempio1}, the vectors $e_i+\eta_i$ span a subspace $W$ of $\ovV$. The corresponding subspace $[W]$ of $\PG(\ovV)$ is a generator of $\cS_q$ opposite to $[V]$ but it is not a complement of $[V]$.  
\end{ex}

\subsection{Regular hermitian varieties}\label{Ex her}

Assume that $\FF$ admits a separable quadratic extension and let $\sigma$ be an involutory non-trivial automorphism of $\FF$. Assuming that $V$ is a right $\FF$-vector space, its dual $V^*$ should be regarded as a left vector space. However, we can re-define it as a right vector space by a vector-times-scalar multiplication defined as follows: $\alpha\cdot t := t^\sigma \alpha$ for every $\alpha\in V^*$. Accordingly, $(a\oplus \alpha)t = at\oplus \alpha t = at\oplus t^\sigma\alpha$. Put $\ovV = V\otimes V^*$, with $V^*$ re-defined as said above. We can define a non-degenerate $\sigma$-hermitian form $h:\ovV\times\ovV\rightarrow \FF$ by declaring that 
\begin{equation}\label{eq h}
h(a\oplus\alpha, b\oplus\beta) ~ = ~ \alpha(b)+(\beta(a))^\sigma, \hspace{5 mm} \forall a, b \in V, ~\alpha, \beta\in V^*.
\end{equation}
The form $h$ is non-degenerate. Let $\cS_h$ be the polar space associated to it. The inclusion mapping of $\cS_h$ in $\PG(\ovV)$ is the unique embedding of $\cS_h$. The subspaces $[V]$ and $[V^*]$ are opposite generators of $\cS_h$ and $\langle [V], [V^*]\rangle = \cS_h$. Hence the unique embedding of $\cS_h$ is tight. 

The hyperbolic lines of $\cS_h$ are the intersections $\ell\cap\cS_h$ for $\ell$ a projective line of $\PG(\ovV)$ which meets $\cS_h$ in at least two points but is not fully contained in $\cS_h$. For instance, if $a\in V$ and $\alpha\in V^*$ are such that $\alpha(a) = 1$ then the hyperbolic line through $[a]$ and $[\alpha]$ consists of $[a]$ and all points $[at\oplus\alpha]$ for $t+t^\sigma= 0$. If a plane $P$ of  $\PG(\ovV)$ contains a point $p\in \cS_h$ such that $p^{\perp_h} \supseteq P$, then either $P\subseteq \cS_h$ of $P\cap\cS_h$ is the union of a set of lines through $p$ and $\ell\cap p^{\perp_h}$ is a hyperbolic line for every line $\ell$ of $P$ which does not pass through $p$. By these fact we can prove an analogue of Proposition \ref{Q2}: if $N$ is a non-deep sub-generator of $\cS_h$ then $N^\perp$ is a hyperbolic line. Consequently, $\cS_h$ is regular. 

Non-deep subgenerators also exist in $\cS_h$, as in $\cS_f$ and $\cS_q$.  The results of Section \ref{Ex symp gen} remain valid for $\cS_h$, with obvious minor modifications. For instance, every generator of $\cS_h$ admits a complement but not all pairs of opposite generators are complementary. Non complementary pairs of opposite generators can be obtained by a minor modifcation of the construction described in Example \ref{esempio2}. 

\subsection{A family of non-regular quadrics}\label{Ex more}

The construction of Section \ref{Ex symp} can be generalized by replacing $V^*$ with a subspace $V'$ of $V^*$ such that $\ker(V') = 0$, possibly $V' \cong V$. We still obtain a symplectic space. As all symplectic spaces are regular \cite{CCGP}, all polar spaces obtained in this way are still regular. In contrast, if we do the same in Sections \ref{Ex hyp} and \ref{Ex her} then the polar space we obtain admits a tight embedding, since it lives in $\PG(V\oplus V')$ and admits $[V]$ and $[V']$ as generators, but in general it is not regular, as in the following case.  

Let $\chr(\FF) = 2$. Given a basis $(e_i)_{i\in\mathfrak{I}}$ of $V$, define the linear functionals $\eta_i$ as in Example \ref{esempio1}. So, $\eta_i(e_j) = \delta_{i,j}$ for every choice of $i,j\in \mathfrak{I}$. Put $V' := \langle \eta_i\rangle_{i\in \mathfrak{I}}$ and $\ovV' := V\oplus V'$. Clearly, $\ker(V') = \{0\}$. With $q$ and $f_q$ as in Section \ref{Ex hyp}, let $q'$ and $f'$ be their restrictions to $\ovV'$ and $\ovV'\times \ovV'$. Then $f'$ is the sesquilinearized of $q'$. 

Explicitly, the vectors of $V\oplus V'$ are sums $v = \sum_{i\in I}e_it_i \oplus \sum_{j\in J}s_j\eta_j$ for $I$ and $J$ finite (possibly empty) subsets of $\mathfrak{I}$, with $t_i \neq 0\neq s_j$ for every $i\in I$ and $j\in J$, by convention. With this convention, we call the ordered  pair $(I,J)$ the {\em support} of $v$. With this convention, if $v = \sum_{i\in I}e_it_i \oplus \sum_{j\in J}s_j\eta_j$ and $v' = \sum_{i\in I'}e_it'_i\oplus\sum_{j\in J'}s'_j\eta_j$ with supports $(I,J)$ and $(I',J')$ respectively. then
\[f'(v,v') = \sum_{i\in I\cap J'}t_is'_i + \sum_{j\in J\cap I'}s_jt'_j ~\mbox{ and } ~ q'(v) = \sum_{i\in I\cap J}t_is_i.\]
Clearly, $f(v,v') = 0$ for every $v$ if and only if $v' = 0$, namely $I' = J' = \emptyset$. Therefore $f'$ is non-degenerate. Hence $q'$ is non-degenerate as well. 

Let $\cS_{q'}$ be the polar space associated to $q'$. As $f'$ is non-degenerate, the inclusion mapping of $\cS_{q'}$ in $\PG(\ovV')$ is the unique embedding of $\cS_{q'}$. Clearly, $[V]$ and $[V']$ are opposite generators of $\cS_{q'}$ and their join spans $\PG(\ovV')$. So, the natural embedding of $\cS_{q'}$ in $\PG(\ovV')$ is tight. 

\subsubsection{Deep and hyperbolic sub-generators} 

We have $\dim(\ker^*(H)\cap V') \leq 1$ for every hyperplane $H$ of $V$. Therefore, if $H$ is a hyperplane of $V$ and $[a\oplus\alpha]\in [H]^\perp$ then either $\alpha = 0$ (hence $a\oplus\alpha \in V$) or $\ker(\alpha) = H$. In the latter case, the condition $\alpha(a) = 0$ forces $a\in H$. Hence $\langle H, a\oplus \alpha\rangle = \langle H, \alpha\rangle = \langle H, \ker^*(H)\rangle$. It follows that $[H]$ is contained in at most two generators of $\cS_{q'}$, namely $[V]$ and possibly $[H+(\ker^*(H)\cap V')]$ (if $\ker^*(H)\cap V' \neq \{0\}$). 

Similarly, since $\dim(\ker(H')) \leq 1$ for every hyperplane $H'$ of $V'$, we obtain that every hyperplane of $[V']$ is contained in at most two generators of $\cS_{q'}$. So, both $[V]$ and $[V']$ are hyperbolic generators. Since every hyperbolic generator contains hyperbolic sub-generators, $\cS_{q'}$ admits hyperbolic sub-generators. Hyperbolic sub-generators are $\perp$-minimal and, if $N$ is a hyperbolic sub-generator and $a, b$ are opposite points of $N^\perp$, then $\{a,b\}^{\perp\perp} = \{a,b\}$. 

So, $\cS_{q'}$ admits hyperbolic lines of size 2. Hence all hyperbolic lines of $\cS_{q'}$ have size 2, since $\Aut(\cS_{q'})$ acts transitively on the set of pairs of opposite points of $\cS_{q'}$. However, as we shall prove in the next subsection (Corollary \ref{Q'3}), $\cS_{q'}$ also admits sub-generators which are neither deep nor hyperbolic. They cannot be $\perp$-minimal. Therefore,

\begin{prop}\label{Q'}
The polar space $\cS_{q'}$ is not regular.
\end{prop} 

\subsubsection{Sub-generators which are not $\perp$-minimal}\label{Ex more 1}  

For $i\in\mathfrak{I}$ put $u_i := e_i\oplus\eta_i$ and, for $i, j\in \mathfrak{I}$, put $u_{i,j} := u_i+u_j$. So, $U := \langle u_{i,j}\rangle_{i,j\in\mathfrak{I}}$ is a maximal totally $q'$-singular subspace of $\ovV'$. Note that, however, $U$ is not maximal among the totally $f'$-isotropic subspaces of $\ovV'$. Indeed $U$ is a hyperplane of $U' := \langle u_i\rangle_{i\in\mathfrak{I}}$ ($= \langle U, u_k\rangle$ for any $k\in\mathfrak{I}$), which is a maximal totally $f$-isotropic subspace of $\ovV'$ (but not a totally $q'$-singular subspace, since none of the vectors $u_i$ is $q'$-singular).   

Pick two elements $i_0$ and $i_1$ of $\mathfrak{I}$, henceforth called $0$ and $1$ for short, and let $\leq$ be a well ordering of $\mathfrak{I}$ with $0$ and $1$ as the first and second element. Put $U_{0,1} := \langle u_{i,j}\rangle_{1 < i < j}$. So, $\cod_U(U_{0,1}) = 2$, namely $[U_{0,1}]$ is a sub-sub-generator of $\cS_{q'}$. Indeed $U = \langle U_{0,1}, u_{k,1}, u_{h,0}\rangle$ for any choice of $h,k\in \mathfrak{I}$ such that $k\neq 1$, $h\neq 0$ and $\{k,h\} \neq \{0,1\}$.  

\begin{lemma}\label{Q'1}
We have $U_{0,1}^{\perp_{f'}} = \langle U_{0,1}, e_0, e_1, \eta_0, \eta_1, u_k\rangle$, for any choice of $k > 1$. 
\end{lemma}
{\bf Proof.} Pick $k \in \mathfrak{I}\setminus\{0,1\}$. Clearly $U_{0,1}^{\perp_{f'}} \supseteq \langle E_{0,1}, e_0, e_1, \eta_0, \eta_1, u_k\rangle$. 
Conversely, let $v = \sum_{i\in I}e_it_i \oplus \sum_{j\in J}s_j\eta_j$ with support $(I,J)$ and suppose that $v$ belongs to $U_{0,1}^{\perp_{f'}}$. We shall prove that $v\in\langle U_{0,1}, e_0, e_1, \eta_0, \eta_1, e_k\oplus\eta_k\rangle$.

Up to adding a suitable combination of $e_0, e_1, \eta_0$ and $\eta_1$ we can assume that $I, J \subseteq \mathfrak{I}\setminus\{0,1\}$. So, let $I, J \subseteq \mathfrak{I}\setminus\{0,1\}$. If $I \neq J$ then we can always find $i, j > 1$ such that $j\not\in I\cup J$ and $i$ belongs to only one of $I$ or $J$.  With $i$ and $j$ chosen in this way, let $i\in I\setminus J$, to fix ideas. Then $f'( u_{i,j}, v)  = t_i$. However $t_i\neq 0$ since $i\in I$ and $(I,J)$ is the support of $v$. It follows that $v\not\perp_{f'}u_{i,j} \in E_{0,1}$. This contradicts the assumption that $v\in E_{0,1}^{\perp_{f'}}$.   

Therefore $I = J$. Suppose now that, for at least one $i \in I$, we have $t_i \neq s_i$. Then, with $j\in \mathfrak{I}\setminus(\{0,1\}\cup I)$, we get  $f'( u_{i,j}, v)  = t_i+s_i \neq 0$. Again, this contradicts the assumption $v\in E_{0,1}^{\perp_{f'}}$. Therefore $t_i = s_i$ for every $i \in I$, namely $v = \sum_{i\in I}u_it_i \in U'$. 

We shall now prove that $v\in \langle U_{0,1}, u_k\rangle$. We shall argue by induction on $|I|$. If $I = \emptyset$ there is nothing to prove. Let $|I| = 1$, say $I = \{i\}$. So, $v = u_it_i$. If $i = k$ there is nothing to prove. Let $i\neq k$. Then $v + u_{i,k}t_i = u_kt_i$. Hence $v = u_kt_i + u_{k,h}t_i \in \langle U_{0,1}, u_k\rangle$, as claimed.

Assuming that $|I | > 1$, let $i$ and $j$ be distinct elements of $I$. Then $v + u_{i,j}t_i$ has support $(J,J)$ with $J$ equal to either $I\setminus\{i,j\}$ or $I\setminus\{i\}$ according to whether $t_i = t_j$ or $t_i \neq t_j$. By the inductive hypothesis, $v+ u_{i,j}t_i \in \langle U_{0,1}, u_k\rangle$. Hence $v\in \langle U_{0,1}, u_k\rangle$. \hfill $\Box$ 

\begin{prop}\label{Q'2}
The star $\St([U_{0,1}])$ of $[U_{0,1}]$ in $\cS_{q'}$ is isomorphic to the quadric $\cQ_4(\FF)$ of $\PG(4,\FF)$.
\end{prop} 
{\bf Proof.} The star $\St([U_{0,1}])$ is isomorphic to the polar space associated with the form induced by $q'$ on a complement $X$ of $U_{0,1}$ in $U_{0,1}^{\perp_{f'}}$. By Lemma \ref{Q'1}, we can choose $X = \langle e_0, e_1, \eta_0, \eta_1, u_k\rangle$. It is straightforward to check that, with $X$ chosen in this way, the restriction of $q'$ to $X$ defines a copy of $\cQ_4(\FF)$. Note that $[u_k]$ is the nucleus of this quadric.  \hfill $\Box$ 

\begin{cor}\label{Q'3}
Each of the sub-generators of $\cS_{q'}$ containing $[U_{0,1}]$ is contained in exactly $|\FF|+1$ (hence at least three) generators of $\cS_{q'}$. 
\end{cor}

For every totally $q'$-singular subpace $W$ of $\ovV'$ containing $U_{0,1}$, the generator $[W]$ of $\cS_{q'}$ is a sub-generator of the polar space $\cS_{f'}$ associated to $f'$. Explicitly, if $W' := \langle W, u_k\rangle$ then $[W']$ is the (unique) generator of $\cS_{f'}$ which contains $[W]$.  (Recall that $[U_{0,1}+\langle u_k\rangle]$ is the nucleus of the quadric $\St([U_{0,1}])$.) 

\begin{prop}
None of the generators of $\cS_{q'}$ containing $[U_{0,1}]$ admits a complement in $\cS_{q'}$. 
\end{prop}
{\bf Proof.} Let $X$ be a complement of $W$ in $\ovV'$, with $W\supseteq  U_{0,1}$ and $[W]$ a generator of $\cS_{q'}$ and let $W' = \langle W, u_k\rangle$. So, $[W']$ is the generator of $\cS_{f'}$ which contains $[W]$. None of the vectors of $W'\setminus W$ is singular for $q'$. Morever, $u_k = w+x$ for suitable vectors $w\in W$ and $x\in X$. So, $x = u_k-w \in W'\setminus W$. Hence $q'(x) \neq 0$. Therefore $[X]$ is not a singular subspace of $\cS_{q'}$. \hfill $\Box$  

\subsubsection{Stars of sub-sub-generators}

Let $X$ be a sub-sub-generator of $\cS_{q'}$. In principle, $\St(X)$ could be any of the following:
\begin{itemize}
\item[(1)] a copy of the hyperbolic quadric $\cQ^+_4(\FF)$ of $\PG(3,\FF)$;
\item[(2)] a copy of the quadric $\cQ_4(\FF)$ of $\PG(4,\FF)$;
\item[(3)] a non-degenerate quadric of rank 2 in a projective space of dimension at least 5 over $\FF$;  
\item[(4)] a quadratic cone in a projective space of dimension at least 3, with a non-degenerate quadric of rank $1$ as a basis;
\item[(5)] a pair of lines meeting at a point;
\item[(6)] one single line. 
\end{itemize}
In case (1) all sub-generators containing $X$ are hyperbolic while in cases (2) and (3) none of them is hyperbolic nor deep. In case (4) the sub-generator corresponding to the vertex of the cone is neither hyperbolic nor deep while all remaining points of $\St(X)$ are provided by deep sub-generators. In case (5) the meet-point of the two lines of $\St(X)$ corresponds to a hyperbolic sub-generator of $\cS_{q'}$ and the remaining points of $\St(X)$ correspond to deep sub-generators. Finally, in case (6) all sub-generators containing $X$ are deep. 

Case (1), (2), (5) and (6) actually occur in $\cS_{q'}$. We have shown in Section \ref{Ex more 1} how to define a sub-sub-generator $X$ with $\St(X)$ as in case (2). In order to obtain (1), choose two hyperplanes $H$ and $H'$ of $V$ (or of $V'$) such that $\ker^*(H)\cap V'\neq\{0\}\neq \ker^*(H')\cap V'$ (respectivley $\ker(H) \neq \{0\} \neq \ker(H')$) and put $X = [H\cap H']$. When $\ker^*(H)\cap V' = \{0\} \neq \ker^*(H')\cap V'$ we get (5) and, if $\ker^*(H\cap H')\cap V' = \{0\}$ then we get (6). We do not know if cases (3) and (4) actually occur. 

In cases (1), (5) and (6) all generators of $\cS_{q'}$ containing $X$ are also generators of $\cS_{f'}$ while in cases (2) and (5) they are deep sub-generators of $\cS_{f'}$. Several possibilities for $\St(X)$ are allowed in cases (3) and (4), depending on the dimensions of the projective space $X^{\perp_{f'}}/X$ which hosts the quadric $\St(X)$ and the dimension of the radical of the bilinearization of the quadratic form which defines that quadric. Of course, the range of hypothetical possibilities depends on $\FF$. For instance, if $\FF$ is finite then $\St(X) \cong \cQ_5^-(\FF)$ is the unique possibility in (3) while two subcases are allowed in (4), according to whether the basis of the cone is a conic of $\PG(2,\FF)$ or an ovoid of $\PG(3,\FF)$. In contrast, if $\FF$ is quadratically closed then (3) is impossible and just one subcase survives in (4), where the basis of the cone is a conic. 

If, as we believe, cases (3) and (4) never occur then the generators of $\cS_{q'}$ are partioned in two families, namely the family ${\cal M}_1$ formed by generators of $\cS_{q'}$ which are generators of $\cS_{f'}$ and the family ${\cal M}_2$ formed by the generators of $\cS_{q'}$ which are (deep) sub-generators of $\cS_{f'}$. The family ${\cal M}_1$ contains all hyperbolic generators of $\cS_{q'}$ while ${\cal M}_2$ contains the non-hyperbolic ones. Moreover, a non-hyperbolic generator of $\cS_{q'}$ contains no deep or hyperbolic sub-generator.   
  
\subsubsection{A conjecture}

The quadric $\cS_q$ described in Section \ref{Ex hyp} is regular because it is hyperbolic (Proposition \ref{Q2}). The proof of Proposition \ref{Q2} relies on an analogue of Lemma \ref{Q1}. In the proof of Lemma \ref{Q1} we have exploited the fact that $\ker(\ker^*(X)) = X$ for every subspace $X$ of $V$. If we replace $V^*$ with a subspace $V'$ of $V^*$ then, in order to repeat that proof, we need that $\ker(\ker^*(X)\cap V') = X$ for every subspace $X$ of $V$. This property holds only if $V' = V^*$. This remark suggests the following conjecture. 

\begin{conj}
With $q$ as in Section \ref{Ex hyp}, let $V'$ be any proper subspace of $V^*$ such that $\ker(V') = \{0\}$. Then the polar space associated to the form induced by $q$ on $V\oplus V'$ is non-regular.   
\end{conj}  

Proposition \ref{Q'} testifies in favor of this conjecture.

\end{document}